\documentclass{paper}

\usepackage{amsmath,amsfonts,amsthm,mathtools}
\usepackage[ruled,linesnumbered]{algorithm2e}
\usepackage{pgf,pgfplots}
\usepackage{latexsym}
\usepackage{amssymb} 

\usepackage{booktabs}
\usepackage[T3,T1]{fontenc}
\usepackage{lmodern}
\DeclareSymbolFont{tipa}{T3}{cmr}{m}{n}
\DeclareMathAccent{\sad}{\mathalpha}{tipa}{16}
\newcommand{\happy}[1]{\breve{#1}}

\title{Componentwise accurate Brownian motion computations using Cyclic Reduction}
\author{Giang T. Nguyen, Federico Poloni}

\newtheorem{theorem}{Theorem}
\newtheorem{lemma}[theorem]{Lemma}
\theoremstyle{definition}

\newtheorem{remark}{Remark}
\newtheorem{corollary}[theorem]{Corollary}

\DeclareMathOperator{\diag}{diag}
\DeclareMathOperator{\offdiag}{offdiag}
\DeclareMathOperator{\tridiag}{tridiag}
\newcommand{\m}[1]{\begin{bmatrix}#1\end{bmatrix}}
\DeclarePairedDelimiter{\abs}{\lvert}{\rvert}
\DeclarePairedDelimiter{\norm}{\lVert}{\rVert}
\newcommand{\bs}[1]{\boldsymbol{#1}}
\newcommand{\ones}{\bs{1}}
\newcommand{\zeros}{\bs{0}}

\newcommand{\leqdot}{\mathrel{\dot{\leq}}}
\renewcommand{\mp}{\mathrm{u}}
\usepackage{hyperref}

\begin{document}
\maketitle

\begin{abstract}
Markov-modulated Brownian motion is a popular tool to model continuous-time phenomena in a stochastic context. The main quantity of interest is the invariant density, which satisfies a differential equation associated with the quadratic matrix polynomial $P(z) = Vz^2-Dz +Q$, where the matrices $V$ and $D$ are diagonal and $Q$ is the transition matrix of a discrete-time Markov chain. Its solution is typically constructed by computing an invariant pair of $P(z)$ associated with its eigenvalues in the left half-plane, or by solving the matrix equation $X^2V-XD+Q=0$. We show that these tasks can be solved using a componentwise accurate algorithm based on Cyclic Reduction, generalizing the recently appeared algorithms for the linear case ($V=0$). We give a proof of the numerical stability of our algorithm in the componentwise sense; the same proof applies to Cyclic Reduction in a more general M-matrix setting which appears in other applications such as the modelling of QBD processes.
\end{abstract}

\section{Introduction}
Markov-modulated Brownian motion~\cite{asmussen95,kk95} 
is a popular tool to model continuous-time phenomena in a stochastic context. An MMBM can be described as the pair $\{Y(t), \phi(t)\}_{t \geq 0}$, where $\phi(t)$ is a continuous-time Markov chain on a state space $\mathcal{S} = \{1,2,\dots,n\}$ with rate matrix $Q\in\mathbb{R}^{n\times n}$ ($Q\ones = \zeros$, where $\ones$ and $\zeros$ are the vectors of all ones and zeros, respectively). Whenever $\phi(t)=i \in \mathcal{S}$, $Y(t)$ evolves according to a Brownian motion process with drift~$d_i$ and variance~$\frac{1}{2}v_i \geq 0$.

The main quantity of interest to determine its steady-state behaviour is the invariant density $\bs{p}(x): \mathbb{R}_{\geq 0}\to \mathbb{R}_{\geq 0}^{1\times n}$, which satisfies, with suitable boundary conditions, the differential equation 
\begin{align} \label{ODE}
 \bs{p}''(x) V - \bs{p}'(x) D + \bs{p}(x) Q= \bs{0},
\end{align}
where $V=\diag(v_i)_{i \in \mathcal{S}}$ and $D=\diag(d_i)_{i \in \mathcal{S}}$. The solutions of this ODE are related to the eigenvalues and left eigenvectors of the matrix polynomial
\begin{align} 
	\label{eq:matpol}
 P(z) := Vz^2-Dz +Q.
\end{align}
The solution of probabilistic interest is asymptotically stable, that is, $\bs{p}(x)\to \bs{0}$ when $x\to\infty$. Hence, we are interested in particular in the eigenvalues $\lambda_i$ with $\Re\lambda_i<0$.

Several methods have been suggested in literature to compute this solution. Some are iterative~\cite{LatN} based on Cyclic Reduction; some depend on the eigendecomposition of a linearization~\cite{kk95}, or more generally a block diagonal decomposition~\cite{AgaS}. A few other algorithms rely on finding a special \emph{invariant pair} $(X,U)$, that is, a pair of matrices $X \in \mathbb{R}^{\ell \times \ell}, U \in \mathbb{R}^{\ell \times n}$ satisfying
\begin{align} 
	\label{invpair2}
	X^2UV - XUD + Q = 0.
\end{align}
For instance, Ivanovs \cite[Sect.~3]{Iva10} considers a related problem---determining the steady-state behavior of a two-boundary Markov-modulated Brownian motion process---which can be solved with the same techniques. The author constructs the solution starting from an invariant pair in which
\begin{align} 
	\label{ivaU}
	U = \begin{bmatrix}
	  I & \Psi
	\end{bmatrix},
\end{align}
where the identity block corresponds to the indices $i$ for which $v_{ii}>0$ or $d_{ii}>0$. This invariant pair $(X, U)$ has a probabilistic meaning: $\Psi \geq 0$ is the matrix recording first-return probabilities of the time-reversed process, and $X$ is a subgenerator matrix ($X_{ij}\geq 0$ if $i\neq j$, and $X\ones \leq \zeros$) for the downward-record process.

A special case often considered in literature is when $v_{ii}>0$ for all $i \in \mathcal{S}$. In this case, $U=I$, and \eqref{invpair2} reduces to
\begin{align} 
	\label{eq:eqRhat}
X^2V-XD+Q=0.
\end{align}
This matrix equation has been studied extensively, especially because of its connection to quasi-birth-death processes~\cite{blm05,ram99}.

Another special case interesting in its own is when $V=0$, that is, when the Markov-modulated Brownian motion $\{Y(t),\phi(t)\}$ is a \emph{stochastic fluid model}, also known as a \emph{fluid queue}. The papers~\cite{xxl12} and~\cite{NguP15} deal with this special case, and provide quadratically convergent algorithms for the invariant density, which are \emph{componentwise accurate}. That is, the algorithms can deliver an approximate solution $\widetilde{\bs{p}}$ such that the quantity $\max_{i \in \mathcal{S}} {\abs{\widetilde{\bs{p}}_i-\bs{p}_i}}/{\bs{p}_i}$ is bounded, rather than ${\norm{\bs{p}-\widetilde{\bs{p}}}}/{\norm{\bs{p}}}$. In this informal introduction, $\bs{p}$ refers to the exact value of a vector quantity related to the solution, for instance, the value of $\bs{p}(x)$ at a determined level $x$, and $\widetilde{\bs{p}}$ represents its computed version in machine arithmetic.

Informally, this means that all entries of $\bs{p}$ have the same number of correct significant digits, irrespectively of their magnitude; thus, all components can be computed to a high accuracy. For example, suppose
\begin{align*} 
\bs{p} = \begin{bmatrix}
  1-10^{-15} & 10^{-15}
\end{bmatrix}.
\end{align*} 
In this case, traditional linear algebra algorithms would instead ensure a high accuracy on the large component $\widetilde{\bs{p}}_1$ only, while $\widetilde{\bs{p}}_2$ could vary wildly with few theoretical guarantees: for instance, it could become negative. Componentwise error bounds are particularly meaningful for probability applications, since small components may represent probabilities of catastrophic failure, which have to be assessed carefully.

In this paper, we focus on computing in a componentwise accurate fashion invariant pairs that solve~\eqref{invpair2} (with the additional property~\eqref{ivaU}) and solutions of the matrix equation~\eqref{eq:eqRhat}, which can then be used to compute a solution $\bs{p}(x)$ of~\eqref{ODE}. This problem contains the linear models treated in~\cite{NguP15,xxl12} as a special case ($V=0$); we extend the techniques introduced there and generalize them to the more challenging second-order case.

In particular, the case in which some of the $v_{ii}$ are zero and some are not, requires special attention. To treat it, we use a method related to both the \emph{shift technique}~\cite{HeMeiniRhee,blm05} and the \emph{theory of index reduction} of differential-algebraic equations~\cite{MehrmannKunkelBook}. We use these techniques in a novel way that combines the themes of these two approaches and adds componentwise accuracy and positivity preservation into the picture.

The paper is structured as follows. In Section~\ref{sec:preliminaries}, we introduce most of the concepts needed in the development of the algorithm, including invariant pairs, componentwise accurate algorithms, and Cyclic Reduction. In Section~\ref{sec:discretization}, we present our solution strategy and formulate a solution algorithm, first for the case $\diag(V)>0$ and then in general. In Section~\ref{sec:stability}, we prove the numerical stability of our algorithm in a componentwise sense. Numerical experiments in Section~\ref{sec:numerical} confirm the effectiveness of this approach, and some brief conclusions follow.

To the best of our knowledge, Lemma~\ref{lem:secondcrtriplet} and the fully subtraction-free version of Cyclic Reduction presented in Algorithm~\ref{algo:crt}, together with the proof of its componentwise stability, are new also in the context of discrete-time quasi-birth-death models~\cite{blm05}.

\section{Assumptions and preliminaries} 
\label{sec:preliminaries}
\subsection{Eigenvalues and invariant pairs of matrix polynomials}

For ease of analysis, we make several assumptions to make sure that the problem cannot be simplified further:
\begin{description}
	\item[A1] The matrix $Q$ is irreducible and aperiodic. 
	\item[A2] $V\neq 0$. 
	\item[A3] there is no index $i \in \mathcal{S}$ for which $v_{ii}=d_{ii}=0$. 
\end{description}

Assumption A1 is to eliminate the cases where our problem can be reduced to smaller disjoint cases. If Assumption A2 is not satisfied, methods for the fluid queue case like the one in~\cite{NguP15} can be used. Finally, if Assumption A3 does not hold, we can replace the problem with another one, where such~$i$ is censored out.

Let $\mathbb{C}[z]$ denote the set of polynomials in the variable $z$. We encountered in the introduction the notion of \emph{eigenvalues} and \emph{left eigenvectors} of a degree-$g$ matrix polynomial 	
	\begin{align*} 
	P(z) = P_0 + P_1 z + \dots + P_g z^g\in\mathbb{C}[z]^{n\times n},
	\end{align*} 
that is, scalars $\lambda\in\mathbb{C}$ and row vectors $\bs{u}\in\mathbb{C}^{1\times n}$ such that $\bs{u}P(\lambda)=\zeros$. Under Assumptions A1--A3, $P(z)$ is a \emph{regular} matrix polynomial, that is, the scalar polynomial $\det P(z)\in \mathbb{C}[z]$ is not identically zero, as one can see by considering its highest-degree term; hence its eigenvalues are a well-defined set of $n$ complex numbers counted with multiplicity. When the leading term $P_g$ is singular, $\det P(z)$ has degree strictly lower than $gn$. In this case we say that $\infty$ is an eigenvalue of $P(z)$ with algebraic multiplicity $gn-\deg P(z)$\footnote{Eigenvalues at infinity are a useful algebraic abstraction that makes several counting and transformation arguments work with little or no modification in a more general setting. We do not venture here in the theory of \emph{geometric multiplicity} and \emph{Jordan structure} for matrix polynomials, which is not a trivial task~\cite{GohLR}.}.

An extremely useful tool to deal with multiple eigenvalues simultaneously, both in theory and in numerical practice, is \emph{invariant pairs} \cite{BetK,GohLR,HigK}. For any $\ell\leq gn$, a pair $(X,U)\in\mathbb{C}^{\ell\times \ell} \times \mathbb{C}^{\ell\times n}$ is called a \emph{left invariant pair} for the matrix polynomial $P(z) = P_0 + P_1 z + \dots + P_g z^g\in\mathbb{C}[z]^{n\times n}$ if
\begin{align} 
	\label{geninvpair}
 UP_0 + XUP_1 + X^2UP_2 + \dots + X^g U P_g = 0
\end{align}
and the matrix $\m{U & XU & \cdots & X^{g-1}U}$ has full row rank. It follows from this definition that if $(X,U)$ is a left invariant pair, then so is
\begin{align} 
	\label{eq:qtrans}
(MXM^{-1},MU)\quad \text{for any nonsingular $M \in \mathbb{C}^{\ell\times \ell}$}. 
\end{align}
The reader not acquainted with this concept can consider a simpler case in which $P(z)$ has $gn$ distinct finite eigenvalues; in this case, the invariant pairs for a matrix polynomial are given by
\begin{align*} 
X=\diag(\lambda_1,\lambda_2,\dots,\lambda_\ell), \quad U = \m{\bs{u}_1\\
\bs{u}_2\\ \vdots\\
\bs{u}_\ell},
\end{align*} 
where $\ell \leq gn$ is arbitrary and for each $i=1,2,\dots,\ell$ the row vector $\bs{u}_i$ is a left eigenvector of $P(z)$ with eigenvalue $\lambda_i$, as well as all the pairs obtained from them through the change of basis transformations in~\eqref{eq:qtrans}. Informally speaking, invariant pairs are a tool to deal with several eigenvalues and eigenvectors at the same time without resorting to a Jordan form.

Invariant pairs generalize the concept of solution of polynomial matrix equations: indeed, if $X$ satisfies~\eqref{eq:eqRhat}, then $(X,I)$ is a left invariant pair for $P(z)$. Moreover, the following properties hold~\cite{GohLR}.
\begin{lemma}
 Let $(X,U)$ be a left invariant pair for a regular matrix polynomial $P(z)$. Then, the following properties hold. 
 \begin{enumerate} \label{invpairprops}
  \item The eigenvalues of $X$ are a subset of the finite eigenvalues of $P(z)$ (both counted with their multiplicity);
  \item \label{deflateinvpair} if
  \[
X = \m{X_{11} & 0\\ X_{21} & X_{22}}, \quad U = \m{U_1 \\ U_2},
  \]
with $X_{11}$ and $X_{22}$ square and $U$ partitioned conformably, then $(X_{11},U_1)$ is another invariant pair for $P(z)$;
 \item $(X,U)$ is an invariant pair also for $P(z)S(z)$, where $S(z)\in\mathbb{C}[z]^{n\times n}$ is any other regular matrix polynomial.
 \end{enumerate}
\end{lemma}

\begin{remark}
A feature that distinguishes linear eigenvalue problems (i.e., the case in which $g=1$, or equivalently $V=0$ in the case of~\eqref{eq:matpol}) from polynomial ones is the fact that in a non-linear eigenproblem \emph{eigenvectors do not uniquely determine their associated eigenvalues}. For instance, both $(1,\bs{e}_1)$ and $(2,\bs{e}_1)$, with $\bs{e}_1=\m{1 & 0}$, are left eigenpairs of the matrix polynomial
\[
P(z) = z^2   I - z \m{3 & 0\\ 0 & 7} + \m{2 & 0\\ 0 & 12}.
\]
Hence, we have to deal explicitly with both elements $U$ and $X$ of the invariant pair, and compute them both at the same time. Instead, in the algorithms for the case $V=0$~\cite{NguP15,xxl12}, it is common to deal with the matrix $\Psi$ in~\eqref{ivaU} as the only unknown, and then compute $X$ from it afterwards. This is possible in the first-order case, but not in the second-order one. This point will prove crucial in Section~\ref{SectionWhereWeNeedToRecoverX}, where the need to compute $X$ as well will impose a nontrivial restriction not present in the linear case.
\end{remark}

\subsection{Triplet representations and accurate matrix exponentials}
As stated earlier, we are interested in performing numerical computations in a way to guarantee the componentwise accuracy of the computed quantities. To this purpose, the main resource are so-called \emph{subtraction-free} algorithms in linear algebra: when the matrices and vectors involved have a prescribed sign structure, it is possible to carry out linear algebraic operations on a computer in terms of sums only, without ever subtracting two floating-point numbers with the same sign. In this case, there is no cancellation, and the results are provably accurate. The most famous algorithm in this class is the \emph{GTH algorithm}~\cite{GraTH85,oci93,AlfXY02} and its generalizations. To introduce it, we need a few preliminary concepts. 

Here and in the following, inequalities between matrices and vectors are used in the componentwise sense: for instance, $A\leq B$ means $A_{ij}\leq B_{ij}$ for each $i,j$.

For a matrix $M\in\mathbb{R}^{n\times n}$, we use the notation $\offdiag(M)$ to denote a vector $\bs{m}\in\mathbb{R}^{n^2-n}$ which contains the elements $\{M_{ij} : i\neq j\}$, i.e., those which do not belong to the main diagonal. The exact ordering of these elements in $\bs{m}$ is not important here. A matrix $M$ is called \emph{M-matrix} if it can be expressed as
\begin{equation} \label{expdecomposition}
M=sI-P, \quad P\in\mathbb{R}_{\geq 0}^{n\times n}, s\geq 0.	
\end{equation}
where $s\in\mathbb{R}$ is greater or equal than the spectral radius $\rho(P)$. It is well known that if an M-matrix $M$ is invertible, then $M^{-1}\geq 0$~\cite{bp94}.

A \emph{triplet representation} for an M-matrix $M$ is a triple $(\bs{m},\bs{v},\bs{w}) \in \mathbb{R}_{\leq 0}^{n^2-n} \times \mathbb{R}_{>0}^{n} \times \mathbb{R}_{\geq 0}^{n}$ such that $\bs{m}=\offdiag(M)$, and $\bs{v}>\zeros$, $\bs{w}\geq \zeros$ are two vectors such that $M\bs{v}=\bs{w}$. The diagonal elements of $M$ do not appear explicitly in the triplet representation, but they are determined uniquely from the relation $M\bs{v}=\bs{w}$. Not all M-matrices admit triplet representations~\cite[Section~1]{guonew}; a counterexample is
$
M = \begin{bsmallmatrix}
  \phantom{-}0 & 0\\
  -1 & 0
\end{bsmallmatrix}.
$
M-matrices that admit a triplet representation are called \emph{regular M-matrices}. 
Non-regular M-matrices must necessarily be singular and reducible~\cite{guonew}, so most M-matrices appearing in applications (and, in particular, all those appearing in the rest of this paper) are indeed regular.

The following result shows that one can solve linear systems with a regular M-matrix with almost perfect componentwise accuracy, given a triplet representation as input.
\begin{theorem}[\cite{oci93,AlfXY02}] 	\label{thm:gth}
Let $(\bs{m},\bs{v},\bs{w})$, be a triplet representation for an invertible regular M-matrix $M\in\mathbb{R}^{n\times n}$, and $\bs{u}\in\mathbb{R}_{\geq 0}^n$. Then, there is an algorithm to compute in $O(n^3)$ floating-point arithmetical operations (starting from the floating-point numbers $\bs{m},\bs{v},\bs{w},\bs{u}$) an approximation $\widetilde{\bs{x}}$ of $\bs{x}=M^{-1}\bs{u} \geq \bs{0}$ such that
\begin{align}
	 \label{gthbound}
\abs{\widetilde{\bs{x}}-\bs{x}} \leq \left(\psi(n)\mp+\mathcal{O}(\mp^2)\right) \bs{x}, 
\end{align}
with $\psi(n)=\frac{2}{3}(2n+5)(n+2)(n+3)$ and $\mp$ being the machine precision.
\end{theorem}

Notice the remarkable absence of the condition number of $M$, which would be necessary in an error bound for an algorithm that uses the matrix entries rather than a triplet representation. The use of a triplet representation as an input makes it possible to solve a linear system with perfect accuracy (up to a polynomial function of the dimension), regardless of its condition number.

The algorithm basically works by computing a LU decomposition of $M$, in which both $L$ and $U$ are M-matrices. Using variants of the same algorithm, one can also perform other related tasks, again starting from a triplet representation $(\bs{m},\bs{v},\bs{w})$ of a regular M-matrix $M$:
\begin{itemize}
	\item computing $M^{-1}$;
	\item solving linear systems of the form $M^\top \bs{x} = \bs{b}$, with $\bs{b}\geq \bs{0}$;
	\item finding the left and right kernel of a singular irreducible $M$.
\end{itemize}

We shall also need the following result.
\begin{lemma} \label{subtriplets}
Let
\[
	\left(\offdiag(M), \begin{bmatrix}
	  \bs{v}_1\\
	  \bs{v}_2
	\end{bmatrix},
	\begin{bmatrix}
	  \bs{w}_1\\
	  \bs{w}_2
	\end{bmatrix}\right), \quad
	M = \begin{bmatrix}
	  M_{11} & M_{12}\\
	  M_{21} & M_{22}
	\end{bmatrix}
\]
(where all the matrices are partitioned conformably) be a triplet representation for the regular M-matrix $M$. Then,
\begin{align}
&(\offdiag(M_{22}), \bs{v}_2, \bs{w}_2 - M_{21}\bs{v}_1), \label{tripletsub}\\
&(\offdiag(S),\bs{v}_1, \bs{w}_1-M_{12}M_{22}^{-1}\bs{w}_2) \label{tripletschur}
\end{align}
are subtraction-free expressions for triplet representations of a principal submatrix $M_{22}$ and its Schur complement $S = M_{11}- M_{12}M_{22}^{-1}M_{21}$.
\end{lemma}
\begin{proof}
The relation $M_{22}\bs{v}_2 = \bs{w}_2 - M_{21}\bs{v}_1$ which defines the first triplet representation comes from expanding the second block row of $M\bs{v}=\bs{w}$. The second relation comes from premultiplying both sides of $M\bs{v}=\bs{w}$ by
\[
\begin{bmatrix}
  I & -M_{12}M_{22}^{-1}\\
  0 & I
\end{bmatrix}.
\]
Additionally, note that $M_{12},M_{21}\leq 0$ and $M_{22}^{-1}\geq 0$ (which can be obtained in a subtraction-free way using the GTH algorithm and the triplet representation~\eqref{tripletsub}), so computing the last terms in~\eqref{tripletsub} and~\eqref{tripletschur} does not involve subtractions. The computation of $S$ via the formula $S=M_{11}- M_{12}M_{22}^{-1}M_{21}$ involves subtractions only for the diagonal entries, but conveniently in the triplet representation we only need $\offdiag(S)$.
\end{proof}
\begin{corollary} 
	\label{diagcorollary}
Given a triplet representation for $M$, $\diag(M)$ can be reconstructed using subtraction-free formulas only.
\end{corollary}
\begin{proof}
It is sufficient to consider~\eqref{tripletsub} in the case in which $M_{22}$ is $1\times 1$. Then, $M_{22} =  (\bs{w}_2 - M_{21}\bs{v}_1) / \bs{v}_2$, where the numerator and denominator are scalar quantities, too.
\end{proof}

We comment briefly also on the computation of the matrix exponential, which we shall need only in the final step of our algorithm. For an M-matrix $M$, it holds that $\exp(-M) \geq 0$. As studied in~\cite{xy08}, it is impossible to find an unconditionally accurate algorithm for this computation in the style of the GTH algorithm; we can only compute approximations $\widetilde{E}$ of $E=\exp(-M)$ satisfying a bound of the form
\[
 \abs{\widetilde{E} - E} = c(M)\mathcal{O}(\mp)E,
\]
where $c(M)$ is a condition number which depends explicitly on $M$. Algorithms for the componentwise accurate computation of matrix exponentials were discussed in~\cite{sgx12}; one of the first steps in these methods is decomposing $M=sI-P$ as in~\eqref{expdecomposition} and using the identity $\exp(-M)=e^{-s}\exp(P)$. Hence it is appealing to look for explicit accurate decompositions of the form~\eqref{expdecomposition} for the matrices whose exponentials we are going to compute.

\subsection{Quadratic matrix equations and their properties}
We now discuss the properties of the solutions of matrix equations of the form $A-BX+CX^2=0$. We focus here on the most common setting in its probabilistic applications; namely, we assume that
\begin{description}
 \item[A4] $A\geq 0$, $C\geq 0$ and $B$ is an M-matrix;
 \item[A5] $(A-B+C)\ones=\zeros$;
 \item[A6] the bi-infinite matrix 
\begin{equation} \label{biinf}
\begin{bmatrix}
  \ddots & \ddots\\
  \ddots & -B & C\\
   & A & -B & C\\
   & & A & -B & \ddots\\
   & & & \ddots & \ddots &
\end{bmatrix}
\end{equation}
is irreducible and aperiodic.
\end{description}
The case most frequently appearing in the probability applications~\cite{blm05} is the one in which $I-B\geq 0$ and $A,I-B,C$ are the transition matrices of a quasi-birth-death (QBD) process, with $A$ being the transition to a lower level.

Under Assumptions A4, A5, A6, one can prove that $A-B+C$ is an irreducible singular M-matrix; we call its left Perron vector $\bs{u}\in\mathbb{R}^{1\times n}$; then we have $\bs{u}>\zeros$, $\bs{u}(A-B+C)=\zeros$. Moreover, one can prove the following results \cite{blm05}.
\begin{theorem} 
	\label{discreteproperties}
Let $A,B,C\in\mathbb{R}^{n\times n}$ satisfying Assumptions A4, A5, A6. Then, the following matrices exist:
\begin{subequations} 
	\label{eq:foursolutions}
  \begin{align}
   G\in\mathbb{R}^{n\times n} &\text{such that $A-BG+CG^2=0$, $G\geq 0$ and $G\ones\leq \ones$};\\
   R\in\mathbb{R}^{n\times n} &\text{such that $R^2A-RB+C=0$, $R\geq 0$ and $\bs{u}R\leq \bs{u}$}.
\end{align}
\end{subequations}
Moreover, the location in the complex plane of the eigenvalues of $F(y)=Ay^2-By+C$, and of $G$ and $R$, is related to the sign of the \emph{mean drift} $d=\bs{u} (C-A) \ones$ as described in Table~\ref{tab:discretecase}.
\begin{table}[h]
\centering
\small
\begin{tabular}{lllllll}
\toprule
Case & Name & $\abs{\mathcal{S}_d}$ & $\abs{\mathcal{U}_d}$ & Other eigvls. & Eigvls. of $G$ & Eigvls. of $R$\\
\midrule
$d<0$ & Positive recurrent & $n$ & $n-1$ & 1 (mult. 1) & $\{\lambda^{-1} : \lambda \in\mathcal{U}_d\} \cup \{1\}$ & $\mathcal{S}_d$  \\
$d=0$ & Null recurrent & $n-1$ & $n-1$ & 1 (mult. 2) & $\{\lambda^{-1} : \lambda \in\mathcal{U}_d\} \cup \{1\}$ & $\mathcal{S}_d \cup \{1\}$  \\
$d>0$ & Transient & $n-1$ & $n$& 1 (mult. 1) & $\{\lambda^{-1} : \lambda \in\mathcal{U}_d\}$ & $\mathcal{S}_d \cup \{1\}$  \\
\bottomrule
\end{tabular}
\caption{Cardinality of the multisets $\mathcal{S}_d = \{\lambda : \lambda$ is an eigenvalue of $F(y)$ and $\abs{\lambda}<1$ $\}$ and $\mathcal{U}_d = \{\lambda : \text{$\lambda$ is an eigenvalue of $F(y)$ and $\abs{\lambda}>1$}\}$, and eigenvalues of the four solution matrices in three possible cases for a triple $A,B,C$ satisfying A4, A5, A6. 
} \label{tab:discretecase}
\end{table}
\end{theorem}
The letters $\mathcal{S}$ and $\mathcal{U}$ in the table stand for `stable' and `unstable', respectively, while $d$ stands for `discrete time' and $c$ in the following will stand for `continuous time'.

Assumption A6 can be relaxed to the less stringent one where $\tridiag(A,-B,C)$ has only one final class \cite[Section~4.7]{blm05}, with only some minor technical complications.

\subsection{Cyclic Reduction}
Cyclic Reduction (CR) is the following matrix iteration. Set 
\begin{align*} 
	A_0 = A, \quad B_0 = \sad{B}_0 = B, \quad C_0 = C,
\end{align*} 
	and compute for each $k=0,1,2,\dots$
\begin{subequations} \label{eq:cr}
  \begin{align}
  A_{k+1} &= A_kB_k^{-1}A_k,\\
  B_{k+1} & = B_k - A_k B_k^{-1} C_k - C_k B_k^{-1} A_k,\\
  C_{k+1} &= C_k B_k^{-1} C_k,\\
  \sad{B}_{k+1} &= \sad{B}_k - C_k B_k^{-1}A_k.
 \end{align}
\end{subequations}

The following applicability and convergence results hold for Cyclic Reduction.
\begin{theorem} \label{thm:crconv}
 Let $A,B,C\in\mathbb{R}^{n\times n}$ satisfying Assumptions A4, A5, A6. Then,
 \begin{enumerate}
  \item $B_k$ is nonsingular for $k\geq 0$; hence, CR can be applied with no breakdown.
  \item $A_k, C_k$ are nonnegative, and $B_k$ and $\sad{B}_k$ are M-matrices for  $k\geq 0$.
  \item $B_k$ and $\sad{B}_k$ converge monotonically to matrices that we shall call $B_\infty$ and $\sad{B}_\infty$, respectively. The matrix $\sad{B}_\infty$ is invertible.
  \item We have
\begin{subequations} 
	\label{eq:crsol}
\begin{align}
 G&=\sad{B}_{\infty}^{-1}A_0,\\
 R&=C_0\sad{B}_{\infty}^{-1}.
\end{align}
\end{subequations}
 \item \label{quadconv} The convergence speed is linear with factor ${1}/{2}$ in the null recurrent case, quadratic with factor $\rho(R)<1$ in the positive recurrent case, and quadratic with factor $\rho(G)<1$ in the transient case.
 \item $(A_k-B_k+C_k)\ones = \zeros$ for each $k\geq 0$, hence $(\offdiag(B_k), \ones, (A_k+C_k)\ones)$ is a triplet representation for $B_k$. \label{crfirsttriplet}
 \end{enumerate}
\end{theorem}
The last item in particular is useful because it allows one to perform the iteration using the GTH algorithm for the inversions required in~\eqref{eq:cr}. Hence, $A_k,\offdiag(B_k),C_k,\offdiag(\sad{B}_k)$ can be computed in a subtraction-free fashion. This is how Cyclic Reduction is currently implemented in software packages such as SMCSolver~\cite{SMCSolver}. However, to implement the final step, \eqref{eq:crsol}, in a subtraction-free way, we need to find a triplet representation for $\sad{B}_\infty$. To this purpose, we give the following result.
\begin{lemma}
	 \label{lem:secondcrtriplet} Under Assumptions A4, A5, A6, the following results hold for the iterates of Cyclic Reduction.
\begin{itemize}
 \item $(A_0 - \sad{B}_k + C_k)\ones = \zeros$ for each $k\geq 0$, hence
 	\begin{align*} 
		\lim_{k\to\infty} C_k\ones=:\sad{\bs{w}} \quad \mbox{ exists, } 	\end{align*} 
		and $(\offdiag(\sad{B}_k),\ones, A_0\ones + \sad{\bs{w}})$ is a triplet representation for $\sad{B}_\infty$.
 \item $\bs{u}(A_k - B_k + C_k)=\zeros$ for each $k\geq 0$.
 \item $\bs{u}(A_k - \sad{B}_k + C_0) = \zeros$ for each $k\geq 0$, hence 
 	\begin{align*} 
		\lim_{k\to\infty} \bs{u}A_k =: \sad{\bs{v}} \mbox{ exists,}
	\end{align*} 
	and $(\offdiag(\sad{B}^{\top}),\bs{u}^{\top},(\bs{u}C_0+\sad{\bs{v}})^{\top})$ is a triplet representation for $\sad{B}_\infty^{\top}$.
\end{itemize}
\end{lemma}
\begin{proof}
We prove only the first equality, the others are analogous. The proof is by induction and similar to the one of item~\ref{crfirsttriplet} of Theorem~\ref{thm:crconv}. For $k=0$, the result holds by Assumption~A5. The inductive step is
 \begin{align*}
  (A_0 - \sad{B}_{k+1} + C_{k+1})\ones &= (A_0 - \sad{B}_{k} + C_kB_k^{-1}A_k +  C_kB_k^{-1}C_k)\ones\\
  &=(A_0 - \sad{B}_{k})\ones + C_kB_k^{-1}(A_k+C_k)\ones\\
  &=(A_0 - \sad{B}_{k})\ones + C_kB_k^{-1}B_k\ones\\
  &=(A_0 - \sad{B}_{k}+C_k)\ones=\zeros,
   \end{align*}
where we have used the fact that $(A_k+C_k)\ones = B_k\ones$ (item~\ref{crfirsttriplet} of Theorem~\ref{thm:crconv}).
\end{proof}
Armed with these triplet representations, we can formulate a fully subtraction-free version of Cyclic Reduction, Algorithm~\ref{algo:crt}.
\begin{algorithm}
\KwIn{$A,B,C\in\mathbb{R}^{n\times n}$ satisfying A4, A5, A6}
\KwOut{The matrices $G$, $R$ defined in~\eqref{eq:foursolutions}.}
Set $A_0=A$, $B_0=\sad{B}_0=B$, $C_0=C$, and $k=0$\;
\Repeat{$\offdiag(\sad{B}_k)$ has converged}{
 Compute $A_{k+1},\offdiag(B_{k+1}), C_{k+1}, \offdiag(\sad{B}_{k+1})$ using \eqref{eq:cr}, performing inversions using the triplet representation in Item~\ref{crfirsttriplet} of Theorem~\ref{thm:crconv}\;
 $k\to k+1$\;
}
Compute $G,R$ using~\eqref{eq:crsol}, performing inversions using one of the triplet representations in Lemma~\ref{lem:secondcrtriplet}\;
\caption{A subtraction-free version of Cyclic Reduction using triplet representations.} \label{algo:crt}
\end{algorithm}

\section{Derivation of the algorithm} \label{sec:discretization}

In this section, we focus on the problem of finding a left invariant pair $(X,U)$ for the matrix polynomial $P(z)$ in~\eqref{eq:matpol} associated to its eigenvalues in the left half-plane. We shall see in Section~\eqref{sec:solvingODE} that a solution to~\eqref{ODE} can be constructed from this pair.

\subsection{The spectrum of $P(z)$}

We start with a theoretical result on the location of the eigenvalues of $P(z)$. To formulate it, we subdivide the indices $i\in\{1,2,\dots,n\}$ into three disjoint subsets, according to the values of $v_{ii}$ and $d_{ii}$, as shown in Table~\ref{tab:states}. Moreover, we set $n_i=\abs{E_i}$ for $i=1,2,3$, so that $n=n_1+n_2+n_3$, and we call $\bs{u}> \zeros$ the left Perron vector of $Q$.

\begin{table}[h]
\centering
\begin{tabular}{lll}
\toprule
Name & $v_{ii}$ & $d_{ii}$\\
\midrule
$E_1$ & $>0$ & any\\
$E_2$ & $=0$ & $>0$\\
$E_3$ & $=0$ & $<0$\\
\bottomrule
\end{tabular}
\caption{Subdivision of each state $i \in \{1,2,\dots,n\}$ into three different sets. Recall that we assume that there is no state with $v_{ii}=d_{ii}=0$.} \label{tab:states}
\end{table}
\begin{theorem} \label{contlocation}
The location in the complex plane of the eigenvalues of $P(z)$ is related to the sign of the \emph{mean drift} $d=\bs{u}D\ones$ as described in Table~\ref{tab:continuouscase}.
\end{theorem}
\begin{proof}
The case $d<0$ appears in~\cite{kk95}; the case $d>0$ can be proved by replacing $D$ with $-D$ (which has the effect of changing the sign of all eigenvalues). 

For the case $d=0$, the proof is not immediate; we give only a sketch, since the result is not necessary for the rest of the paper. A limit argument from both sides shows that $\abs{\mathcal{S}_c}\leq n_1+n_2-1$ and $\abs{\mathcal{S}_u}\leq n_1+n_3-1$. Since we assume that $v_{ii}$ and $d_{ii}$ are not both zero, $D_{E_2 \cup E_3,E_2 \cup E_3}$ is nonsingular, and hence the $n_2+n_3$ Jordan chains for $\lambda=\infty$ (as defined in~\cite[Section~1.4]{GohLR}) have length $1$. Hence the only thing left to prove is that there are no more than $2$ eigenvalues on the imaginary axis. The Gerschgorin argument in~\cite{kk95} shows that the only possible eigenvalue on the imaginary axis is zero. The multiplicity of the eigenvalue $0$ is at most $2$, because for $h$ small enough the matrix polynomial $I+hP(z)$ satisfies the hypotheses of Theorem~\ref{discreteproperties}, as we show in more detail in the following.
 \end{proof}

\begin{table}[h]
\centering
\begin{tabular}{llllll}
\toprule
Case & Name & $\abs{\mathcal{S}_c}$ & $\abs{\mathcal{U}_c}$ & Other eigenvalues\\
\midrule
$d<0$ & Pos. rec. & $n_1+n_2$ & $n_1+n_3-1$ & 0 (mult. 1), $\infty$ (mult. $n_2+n_3$)\\
$d=0$ & Null rec. & $n_1+n_2-1$ & $n_1+n_3-1$ & 0 (mult. 2), $\infty$ (mult. $n_2+n_3$)\\
$d>0$ & Transient & $n_1+n_2-1$ & $n_1+n_3$ & 0 (mult. 1), $\infty$ (mult. $n_2+n_3$)\\
\bottomrule
\end{tabular}
\caption{Cardinality of $\mathcal{S}_c = \{\lambda : \text{$\lambda$ is an eigenvalue of $P(z)$ and $\Re{\lambda}<0$}\}$ and $\mathcal{U}_c = \{\lambda : \text{$\lambda$ is an eigenvalue of $P(z)$ and $\Re{\lambda}>0$}\}$.
} \label{tab:continuouscase}
\end{table}

We now have all we need to define precisely which invariant pair we are looking for. We call a left invariant pair of $P(z)$ \emph{c-stable}, if its associated eigenvalues are:
\begin{itemize}
	\item the $n_1+n_2$ eigenvalues in $\mathcal{S}_c$, if $P(z)$ is positive recurrent; or
	\item the $n_1+n_2-1$ eigenvalues in $\mathcal{S}_c$ and the eigenvalue $0$ with multiplicity $1$, if $P(z)$ is null recurrent or transient.
\end{itemize}
Similarly, with respect to Table~\ref{tab:discretecase}, we call a left invariant pair of $F(y)$ \emph{d-stable}, if its associated eigenvalues are:
\begin{itemize}
	\item the $n$ eigenvalues in $\mathcal{S}_d$, if $F(y)$ is positive recurrent; or
	\item the $n-1$ eigenvalues in $\mathcal{S}_d$ and the eigenvalue $1$ with multiplicity $1$, if $F(y)$ is null recurrent or transient.
\end{itemize}
Notice that $(R,I)$, where $R$ is the matrix in~\eqref{eq:foursolutions}, is a d-stable invariant pair for $F(y)$.

\subsection{The case $\diag(V)>0$} \label{SectionWhereWeNeedToRecoverX}
We start by treating the simpler case in which $\diag(V)>0$ (or, in probabilistic terms, the dynamic in all states has a Brownian motion component). We have $E_2=E_3=\varnothing$ and $n_2=n_3=0$, all $2n$ eigenvalues of $P(z)$ are finite and we are looking for exactly $n$ of them. The formulation in~\eqref{ivaU} has $U=I$, hence the task of finding an invariant pair becomes the one of finding a solution of the matrix equation~\eqref{eq:eqRhat}.

We have seen that Cyclic Reduction can be applied to matrix polynomials $F(y)$ with a specific sign structure, which is associated with a specific spectral structure as shown in Table~\ref{tab:discretecase}. The sign structure and spectral structure of $P(z)$ in~\eqref{tab:continuouscase} do not match these requirements, so we need some preprocessing to convert one case into the other. Even if the sign structure is a stricter requirement, it is useful for our analysis to focus first on the spectral structure, and describe methods of altering the position of the eigenvalues.

We start from a general lemma on rational transformations of matrix polynomials.
\begin{lemma}[\cite{M4poly,Nofpoly}] \label{rattrans}
 Let 
 \begin{align*} 
 	y = f(z)=\frac{\alpha z+\beta}{\gamma z+\delta}
\end{align*} 
be a degree-1 (scalar) rational function, with $\alpha,\beta,\gamma,\delta\in\mathbb{C}$ and $\alpha\delta\neq \beta\gamma$, $z=f^{-1}(y)={(\delta y-\beta)}/{(\alpha-\gamma y)}$ its inverse, and $P(z)\in\mathbb{R}[z]^{n\times n}$ be a degree-$g$ regular matrix polynomial with eigenvalues $\lambda_1,\dots,\lambda_{gn}$ (counted with multiplicity, and possibly including $\infty$). Then, the following properties hold. 
 \begin{enumerate}
  \item The matrix polynomial 
  	\begin{align*} 
		F(y)=(\alpha-\gamma y)^g P({\delta y-\beta}/{\alpha -\gamma y})
	\end{align*} 
	is regular and has eigenvalues $f(\lambda_i)$, for each $i=1,2,\dots,gn$.
  \item If $(Z,U)$ is a left invariant pair for $P(z)$ and $\gamma Z+\delta I$ is nonsingular, then $(f(Z),U)$ is an invariant pair for $F(y)$. Conversely, if $(Y,U)$ is an invariant pair for $F(y)$ and $\alpha I-\gamma Y$ is nonsingular, then $(f^{-1}(Y),U)$ is an invariant pair for $P(z)$.
 \end{enumerate}
\end{lemma}
Note that $f(Z) = (\gamma Z+ \delta I)^{-1}(\alpha Z+\beta I)=(\alpha Z+\beta I)(
\gamma Z+\delta I)^{-1}$ is well-defined for a matrix argument $Z$ since the two factors commute, and similarly for $f^{-1}(Y)$. If $U=I$, the last item gives a relation between the solutions to the unilateral matrix equations associated to $P(z)$ and $F(y)$.

Lemma~\ref{rattrans} suggests a general strategy to approach the problem:
\begin{enumerate}
 \item Choose a function $f$ such that $f(0)=1$ and the images of $\mathcal{S}_c,\mathcal{U}_c$ lie inside and outside the unit circle, respectively.
 \item Construct $F(y)=Ay^2-By+C$ as in Lemma~\ref{rattrans}.
 \item Apply Cyclic Reduction to find the solution $R$ to 
 	\begin{align*} 
		R^2A-RB+C=0;
	\end{align*} 
	then, $(R,I)$ is a d-stable invariant pair of $F(y)$.
 \item Compute $X=f^{-1}(R)$. \label{laststep}
\end{enumerate}
This general framework of relocating eigenvalues via rational transformations is quite common in literature; see for instance~\cite{bmp} for a discussion of it in the case of fluid queues ($V=0$). Frequent choices for $f$ are
	\begin{align*} 
		 y=1+hz, \quad \mbox{and} \quad 
		 	y=\frac{1+hz}{1-hz},
	\end{align*} 
	where $h>0$ is a parameter. However, some care is needed here to allow for componentwise accurate computations within the framework. The first important restriction comes from the last step: once we have obtained $R\geq 0$, we need to be able to compute $f^{-1}(R)$. If one chooses $y={1+hz}/{1-hz}$, the computation becomes $X= h^{-1}(I-R)(I+R)^{-1}$. This is problematic, because the matrix $I+R$, which we need to invert, is a nonnegative matrix; hence Theorem~\ref{thm:gth} does not apply, and we do not know of another componentwise algorithm to invert matrices with this sign pattern, even if triplet representations are available.

Things are easier if one chooses the function $y=1+hz$. In this case, the last step becomes $X= h^{-1}(R- I)$; subtractions are needed only to compute its diagonal, and we can avoid them completely using Corollary~\ref{diagcorollary} if we manage to obtain a triplet representation for $-X$. Moreover, the matrix is now naturally expressed in the form~\eqref{expdecomposition}. For this reason, we set $y=f(z)=1+hz$ in the following.

With this choice, we get $z={(y-1)}/{h}$, and
\begin{align*} 
 P(z) & = V\frac{(y-1)^2}{h^2} -D\frac{y-1}{h} + Q \\
 & = \frac{1}{h^2}Vy^2 - \left(2\frac{1}{h^2}V + \frac{1}{h}D\right)y + \left(\frac{1}{h^2}V+\frac{1}{h}D+Q\right),
\end{align*} 
hence
\begin{align} 
	\label{ABC}
A := \frac{1}{h^2}V, \quad B:= 2\frac{1}{h^2}V + \frac{1}{h}D, \quad C:= \frac{1}{h^2}V+\frac{1}{h}D+Q.	
\end{align}

Once we have decided to use~\eqref{ABC} to convert the setting into that of a discrete-time quadratic matrix equations, we have to choose the value of the parameter~$h$. A first requirement is that Assumption A4 is satisfied; it is easy to see that it holds provided that $v_{ii}-d_{ii}h+q_{ii}h^2 \geq 0$ for each $i$. Since we assume $v_{ii}>0$ for each $i \in \mathcal{S}$ for now, this holds for sufficiently small values of $h$. Moreover, we have to ensure that the computed diagonal of $C$ is componentwise accurate. For this, we follow the strategy used in~\cite{xxl12}: we choose $h$ small enough so that all the required subtractions are of the form $b-a$ with $b\geq 2a\geq 0$. In this case, there cannot be catastrophic cancellation in the subtraction in machine arithmetic. This requirement translates to the following constraints on $h$:
\begin{subequations} \label{eq:hconstraints}
\begin{align}
 v_{ii} &\geq -2(d_{ii}h + q_{ii}h^2), \text{ for each $i$ with $v_{ii}>0$ and $d_{ii}<0$};\\
 v_{ii}+hd_{ii} &\geq -2 q_{ii}h^2, \text{ for each $i$ with $d_{ii}>0$}.
\end{align}
\end{subequations}
All these inequalities are satisfied for a sufficiently small value of $h$, which is easy to compute explicitly. Another possibility is performing these subtractions using machine arithmetic with a higher precision; since there are only $O(n)$ of them, this safeguard will not impact the final cost of the algorithm.

It is easy to check that Assumption A5 is always satisfied. We prove below that A6 is satisfied as well.

\begin{theorem}
Suppose the matrix $Q$ is irreducible and aperiodic, and $V\neq 0$. Then, Assumption A6 holds for the matrices $A,B,$ and $C$ defined in~\eqref{ABC}.
\end{theorem}
\begin{proof}
We identify each element in the index set of the infinite matrix~\eqref{biinf} with a pair $(i,\ell)\in \{1,2,\dots,n\} \times \mathbb{Z}$, where the second entry denotes the block (level) and the first denotes the position in the block. We shall prove that there is a walk in the graph associated to~\eqref{biinf} between any two states $(i,\ell),(j,m) \in \{1,2,\dots,n\} \times \mathbb{Z}$.

Let $k\in\{1,2,\dots,n\}$ be such that $v_{kk}>0$. By the irreducibility assumption, we can find in the graph associated to $Q$ a walk from $i$ to $j$ which passes through $k$ and has length at least $m-\ell$. Since $C$ has the same offdiagonal nonzero structure as $Q$, the same walk can be used in the matrix~\eqref{biinf}, and after each step the second element of the pair goes up by one. Hence the path goes from $(i,\ell)$ to $(j,m+p)$, for some $p>0$. We modify this walk by inserting $p$ transitions using the nonzero entry $A_{kk}$ when we first reach $k$ as the first element of the pair. The resulting graph goes from $(i,\ell)$ to $(j,m)$, as requested.

If $Q$ is aperiodic, the same construction can be made with different values of $p$ which are coprime; hence~\eqref{ABC} is aperiodic, too.
\end{proof}

Moreover, we can prove that our transformation~\eqref{ABC} preserves the sign of the mean drift.
\begin{lemma}
Let $d_c$ be the mean drift of the Markov-modulated Brownian motion process with parameters $V,D,Q$, and $d_d$ be the mean drift of the QBD process associated to $A,B,C$ as in~\eqref{ABC}. Then, $d_d = h^{-1}d_c$.
\end{lemma}
\begin{proof}
First note that $A-B+C=Q$, so the left Perron vector $\bs{u}$ of $A-B+C$ coincides with the one of $Q$. Then it is easy to compute
\[
d_d = \bs{u}(C-A)\bs{1} = \bs{u}\left(\frac{1}{h^2}V+\frac{1}{h}D+Q-\frac{1}{h^2}V\right)\ones = \frac{1}{h}\bs{u}D\ones = {1}{h}d_c.
\]
\end{proof}
Finally, we note that as a byproduct of Cyclic Reduction (Algorithm~\eqref{algo:crt}) we can obtain explicitly a triplet representation for $-X^{\top}$.
\begin{theorem}
The triplet 
\begin{equation} \label{tripletX}
\left(\offdiag(-X^{\top}),\bs{u}^{\top},\left(h^{-1}\sad{\bs{v}}\sad{B}^{-1}_\infty\right)^{\top}\right)
\end{equation}
is a triplet representation for the matrix $-X^{\top}$, where $X=h^{-1}(R-I)$.
\end{theorem}
\begin{proof}
By Lemma~\ref{lem:secondcrtriplet}, we have $\sad{\bs{v}} - \bs{u}\sad{B}_\infty + \bs{u}C_0 = \zeros$. Hence, 
	\begin{align*} 
		-\bs{u}X = \bs{u}\frac{1}{h}(I-C_0\sad{B}_{\infty}^{-1}) = \frac{1}{h}\sad{\bs{v}}\sad{B}^{-1}_\infty.
	\end{align*}
\end{proof}

Summing up everything, our algorithm for the case $\diag(V)>0$ is described in Algorithm~\ref{algo:vpos}.
\begin{algorithm}
\KwIn{$V,D\in\mathbb{R}^{n\times n}$ diagonal matrices with $\diag(V)>\zeros$, $Q\in\mathbb{R}^{n\times n}$ a generator matrix ($Q\ones =\zeros$, $\offdiag(Q)\geq \zeros$), satisfying A1, A2, A3.}
\KwOut{the matrix $X$ (or a decomposition~\eqref{expdecomposition} for it), and the triplet representation~\eqref{tripletX} for $-X^{\top}$.}
compute the left Perron vector $\bs{u}$ of $Q$ using the triplet representation $(\offdiag(-Q),\ones,\zeros)$\;
choose $h$ small enough so that $h^{-2}v_{ii}+h^{-1}d_{ii}+q_{ii}>0$, and it can be computed without catastrophic cancellation\;
compute $A,B,C$ as in~\eqref{ABC}\;
compute $R \geq 0$ via Algorithm~\ref{algo:crt}\;
using the last iterate $A_k$ computed by Algorithm~\ref{algo:crt} and the triplet representation for $\sad{B}_{\infty}$, compute $\sad{\bs{v}}=\bs{u}A_k$ and $h^{-1}\sad{\bs{v}}\sad{B}_{\infty}^{-1}$\;
compute $X=h^{-1}(R-I)$ (or $P= h^{-1}R$ and $s = h^{-1}$)\;
\caption{Computing a c-stable invariant pair $(X,I)$ of $P(z)$, in the case $\diag(V)>0$} \label{algo:vpos}
\end{algorithm}

\subsection{Shifting infinite eigenvalues in $P(z)$}
The method outlined in the previous section uses the assumption that $v_{ii}>0$ for each $i$. When this is not the case, it is not true in general that we can choose~$h$ small enough to have $h^{-2}v_{ii}+h^{-1}d_{ii}+q_{ii}\geq 0$. This is possible for $i\in E_2$, since $d_{ii}>0$, but if $E_3$ is not empty the algorithm cannot be applied. 

Moreover, if $V$ is singular, then the matrix polynomial $P(z)$ has infinite eigenvalues, and a c-stable invariant pair $(X,I)$, with $X$ of size $n\times n$, cannot be constructed since even in the positive recurrent case $P(z)$ does not have $n$ eigenvalues in the left half-plane. 

Finally, the discretization methods outlined in the previous section all break down in some way: if we use the map $y=1+hz$, then we cannot enforce the requirement that the eigenvalue $z=\infty$ is mapped inside the unit circle by choosing a small enough $h$; if we use a variant of the Cayley transform, then $f(\infty)=-1$, and we are left with an eigenvalue of $F(y)$ at $-1$, possibly with high multiplicity; this eigenvalue often prevents the convergence of Cyclic Reduction (note indeed that we are not in the hypotheses of Theorem~\ref{thm:crconv}). All these issues are related, and indeed we can solve all of them with the same modification to the algorithm.

We subdivide the parameter matrices into blocks corresponding to $E_1,E_2,E_3$ as
\begin{equation} \label{VDQ3blocks}
 V = \m{V_1 & 0 & 0\\ 0 & 0 & 0\\0&0&0}, \quad D = \m{D_1 & 0 & 0\\ 0 & D_2 & 0\\0&0& D_3}, \quad Q  = \m{Q_{11} & Q_{12} & Q_{13}\\Q_{21} & Q_{22} & Q_{23}\\Q_{31} & Q_{32} & Q_{33}},	
\end{equation}
where $V_1>0$, $D_2 > 0$ and $D_3 \leq 0$ are diagonal matrices.

We define
\[
\widetilde{P}(z) := P(z)S(z), \quad S(z) := \m{I & 0 & 0\\ 0 & I & 0\\ 0 & 0 & (1+hz)I}.
\]
The resulting matrix polynomial $\widetilde{P}(z)=\widetilde{V}z^2-\widetilde{D}z+\widetilde{Q}$ has coefficients
\[
 \widetilde{V} := \m{V_1 & 0 & 0\\ 0 & 0 & 0\\0&0&-hD_{3}}, \quad \widetilde{D} := \m{D_1 & 0 & -hQ_{13}\\ 0 & D_2 & -hQ_{23}\\0&0& D_3-hQ_{33}}, \quad \widetilde{Q} := Q.
\]
Every finite eigenvalue $\lambda\neq\infty$ of $P(z)$ is also an eigenvalue of $\widetilde{P}(z)$ (with the same left eigenvector), while $n_3$ infinite eigenvalues are replaced by eigenvalues $-h^{-1}$. This can be readily proved by considering the determinants $\det \widetilde{P}(z) = \det P(z) \det S(z)$ and their degrees. This formulation of shifting as multiplication by a suitable matrix polynomial has been suggested recently in~\cite{BinM_brauer}.

\begin{remark}
We can interpret this transformation as a manipulation of the differential equation~\eqref{ODE}. Indeed, if we subdivide $\bs{p}(x) = \begin{bmatrix}
  \bs{p}_1(x) & \bs{p}_2(x) & \bs{p}_3(x)
\end{bmatrix}$ 
conformably, then the third block equation reads 
\begin{align} 
	\label{ODE3}
 - D_3\bs{p}_3'(x) + Q_{31}\bs{p}_1(x)+Q_{32}\bs{p}_2(x)+Q_{33}\bs{p}_3(x) = \bs{0};
\end{align}
differentiating this equation gives
\begin{align} 
	\label{ODE3d}
- D_3\bs{p}_3''(x) + Q_{31}\bs{p}_1'(x)+Q_{32}\bs{p}_2'(x)+Q_{33}\bs{p}_3'(x) = \bs{0};
\end{align}
then the equation $\bs{p}''(x)\widetilde{V} - \bs{p}'(x)\widetilde{D}+\bs{p}(x)Q = \bs{0}$ is obtained from~\eqref{ODE} by replacing the third block equation~\eqref{ODE3} with $\eqref{ODE3} + h\eqref{ODE3d}$. This kind of manipulations is commonly used in the context of index reduction techniques~\cite{MehrmannKunkelBook}.
\end{remark}
We can set up the discretization scheme described in Section~\ref{sec:discretization} starting from $\widetilde{P}(z)$ rather than $P(z)$. The resulting polynomial $\widetilde{F}(y)$ has coefficients
\begin{subequations} \label{eq:tildeABC}
\begin{align}
\widetilde{A} &:= \m{h^{-2}V_1 & 0 & 0\\ 
0 & 0 & 0\\
0&0&- h^{-1}D_{3}},\\
\widetilde{B}&:=\m{ 2h^{-2}V_1 + h^{-1}D_1 & 0 & -Q_{13}\\ 0 & h^{-1}D_2 & -Q_{23}\\0&0& -h^{-1}D_3-Q_{33}},\\
\widetilde{C}&:= \m{h^{-2}V_1+ h^{-1}D_1+Q_{11} & Q_{12} & 0\\ Q_{21} & h^{-1}D_2 + Q_{22} & 0\\Q_{31}&Q_{32}&0}.
\end{align}
\end{subequations}
Notice the nontrivial simplification that zeroes out the last block column of $\widetilde{C}$. Its appearance is due to the fact that the eigenvalues $-h^{-1}$ introduced in the previous step get mapped to $f(-h^{-1})=0$, hence $\widetilde{F}(y)$ has $n_3$ zero eigenvalues.
 
If one chooses a sufficiently small $h$, the diagonals of $h^{-2}V_1+ h^{-1} D_1+Q_{11}$ and $h^{-1}D_2 + Q_{22}$ are nonnegative (and can be computed accurately): it is sufficient to impose~\eqref{eq:hconstraints} on $i\in E_1 \cup E_2$. Hence, the matrices $\widetilde{A}$ and $\widetilde{C}$ are nonnegative, and $\widetilde{B}$ is an M-matrix, which is the correct sign structure to implement subtraction-free Cyclic Reduction (Algorithm~\ref{algo:crt}).

\subsection{Deflating zero eigenvalues in $R$}  \label{sec:deflatingzeros}
In the case $v_{ii}$ may be zero, the solution $R$ produced by Cyclic Reduction does not give immediately the invariant pair we need. Indeed, in view of our previous analysis of the eigenvalues of $F(y)$ and $\widetilde{F}(y)$, in the positive recurrent case the eigenvalues of $R$ comprise of 
\begin{itemize}
 \item $n_3$ zero eigenvalues;
 \item $f(\lambda)$, for each eigenvalue $\lambda$ of $P(z)$ with $\Re \lambda<0$, counted with multiplicity.
\end{itemize}
The eigenvalues of the form $f(\lambda)$ are precisely the ones  we need in our invariant pair, but there are spurious zero eigenvalues. If $R$ were in the form
\begin{equation} \label{trailingzerocolumns}
\begin{bmatrix}
  * & 0\\
  * & 0
\end{bmatrix},
\end{equation}
with the bottom-right block $n_3\times n_3$, we could remove them by applying the result in point~\ref{deflateinvpair} of Theorem~\ref{invpairprops} to the invariant pair $(R,I)$.

Unfortunately, we have $R=C_0 \sad{B}_\infty^{-1}$, where $C_0$ is in the form~\eqref{trailingzerocolumns} and $\sad{B}_\infty$ is a regular M-matrix (for which we know a triplet representation). When one carries out the product, the zero block is lost. To recover it, we have to switch to a different invariant pair.
\begin{theorem}
Let the matrix $C$ in~\eqref{eq:tildeABC} and the matrix $\sad{B}_\infty$ produced by Cyclic Reduction on~\eqref{eq:tildeABC} be partitioned as
\begin{align} 
	\label{BCpartitioning}
C_0 = \begin{bmatrix}
  C_{11} & 0\\
  C_{21} & 0
\end{bmatrix}, \quad
\sad{B}_\infty =\begin{bmatrix}
  B_{11} & B_{12}\\
  B_{21} & B_{22}
\end{bmatrix},	
\end{align}
where the bottom-right block of dimension $n_3\times n_3$ corresponds to the indices in $E_3$. Then, $(Y,\begin{bmatrix}
	    I & \Psi
	\end{bmatrix})$, with
\begin{subequations} \label{dstableinvpairformula}
\begin{align*}
	\Psi & = -B_{12}B_{22}^{-1} \geq 0,\\
	Y & = (C_{11}+\Psi C_{21})S^{-1}\geq 0, \quad  \text{with } S := B_{11}+\Psi B_{21}
\end{align*}
\end{subequations}
is a subtraction-free expression for a d-stable left invariant pair of $\widetilde{F}(y)$, and 
\begin{equation} \label{cstableinvpairformula}
(X,\begin{bmatrix}
	    I & \Psi
	\end{bmatrix}), \quad \text{with } X= h^{-1}(Y-I)
\end{equation}
is a subtraction-free expression for a c-stable left invariant pair of $P(z)$.
\end{theorem}
\begin{proof}
We apply to the d-stable left invariant pair $(R,I)$ of $\widetilde{P}(z)$ a transformation of the form~\eqref{eq:qtrans} with
\[
M = \begin{bmatrix}
  I & \Psi\\
  0 & I
\end{bmatrix},
\]
obtaining
\begin{align*}
(MRM^{-1},M) &= (MC(M\sad{B}_{\infty})^{-1},M)\\
& = \left(
\begin{bmatrix}
  C_{11}+\Psi C_{21} & 0\\
  C_{21} & 0
\end{bmatrix}
\begin{bmatrix}
  B_{11}+ \Psi B_{21} & 0\\
  B_{21} & B_{22}
\end{bmatrix}^{-1},
\begin{bmatrix}
  I & \Psi\\
  0 & I
\end{bmatrix}
\right)\\
& = \left(
\begin{bmatrix}
  (C_{11}+\Psi C_{21})S^{-1} & 0\\
  C_{21}S^{-1} & 0
\end{bmatrix},
\begin{bmatrix}
  I & \Psi\\
  0 & I
\end{bmatrix}
\right).\\
\end{align*}
Notice that $B_{22}$ and $S$ are respectively a submatrix and a Schur complement of the regular M-matrix $\sad{B}_\infty$, so triplet representations to invert both are available by Lemma~\ref{subtriplets}. We can now apply point~\ref{deflateinvpair} of Theorem~\ref{invpairprops} to obtain that~\eqref{dstableinvpairformula} is an invariant pair associated to the d-stable eigenvalues of $F(y)$. Transforming this invariant pair with Lemma~\ref{rattrans}, we obtain~\eqref{cstableinvpairformula}.
\end{proof}

\subsection{A triplet representation for $-X^{\top}$} \label{sec:tripletX2}
In this section, we obtain a triplet representation for the M-matrix $-X^{\top}$ using subtraction-free expressions only.
\begin{theorem}
The triplet
\begin{equation} \label{tripletX2}
	\left(\offdiag(-X^{\top}),\bs{u}_1^{\top}, \frac{1}{h}\left(\left(\sad{\bs{v}}_1 + \sad{\bs{v}}_2B_{22}^{-1}(C_{21}-B_{21})\right)S^{-1}\right)^{\top}\right)
\end{equation}
is a subtraction-free expression for a triplet representation of $-X^{\top}$, where $X$ is defined by~\eqref{cstableinvpairformula}.
\end{theorem}
\begin{proof}
By introducing the partitioning~\eqref{BCpartitioning} in Lemma~\ref{lem:secondcrtriplet}, we get $\bs{u}_1B_{12} + \bs{u}_2 B_{22} = \sad{\bs{v}}_2$. Hence, $\bs{u}_2 = \sad{\bs{v}}_2 B_{22}^{-1} + \bs{u}_1\Psi $, and
\begin{equation} \label{vBeq}
	\bs{u}M^{-1} = \begin{bmatrix}
	  \bs{u}_1 & \bs{u}_2
	\end{bmatrix}
	\begin{bmatrix}
	  I & -\Psi\\
	  0 & I
	\end{bmatrix}
	= \begin{bmatrix}
	  \bs{u}_1 & \sad{\bs{v}}_2B_{22}^{-1}
	\end{bmatrix}.
\end{equation}

Again, from Lemma~\ref{lem:secondcrtriplet}, we get
\begin{align*}
\sad{\bs{v}} &= \bs{u}(\sad{B}_\infty -C_0) = \bs{u}M^{-1}(M\sad{B}_\infty -MC_0)\\
&=\begin{bmatrix}
  \bs{u}_1 & \sad{\bs{v}}_2B_{22}^{-1}
\end{bmatrix}
\left(
\begin{bmatrix}
  S & 0\\
  B_{21} & B_{22}
\end{bmatrix}
-
\begin{bmatrix}
	C_{11} + \Psi C_{21} & 0\\
 	C_{21} & 0
\end{bmatrix}
\right),
\end{align*}
and the first block column of this expression gives
\[
\sad{\bs{v}}_1 + \sad{\bs{v}}_2B_{22}^{-1}(C_{21}-B_{21}) = \bs{u}_1 (S - (C_{11}+\Psi C_{21}))
\]
or
\[
\left(\sad{\bs{v}}_1 + \sad{\bs{v}}_2B_{22}^{-1}(C_{21}-B_{21})\right)S^{-1} = \bs{u}_1 (I - (C_{11}+\Psi C_{21})S^{-1}) = -h \bs{u}_1 X,
\]
from which~\eqref{tripletX2} follows. Note that $B_{21}\leq 0$ and $C_{21}, B_{22}^{-1},S^{-1}\geq 0$, so no subtractions are needed in~\eqref{tripletX2}.
\end{proof}
\subsection{The algorithm}

Putting everything together, we obtain Algorithm~\ref{algo:npt} for the computation of the c-stable invariant pair of a matrix polynomial $P(z)$, which generalizes Algorithm~\ref{algo:vpos} by removing the assumption that $\operatorname{diag}(V)>0$.
\begin{algorithm}
\KwIn{$V,D\in\mathbb{R}^{n\times n}$ diagonal matrices with $\diag(V) \geq \zeros$, $Q\in\mathbb{R}^{n\times n}$ a generator matrix ($Q\ones =\zeros$, $\offdiag(Q)\geq \zeros$), satisfying A1, A2, A3.}
\KwOut{a c-stable invariant pair $(X,\begin{bmatrix}
  I & \Psi
\end{bmatrix})$ of $P(z) = Vz^2-Dz+Q$ (or a decomposition~\eqref{expdecomposition} for it) and a triplet representation for $-X^{\top}$.}

compute the left Perron vector $\bs{u}$ of $Q$ using the triplet representation $(\offdiag(-Q),\ones,\zeros)$\;
choose $h > 0$ small enough to satisfy~\eqref{eq:hconstraints}\;
compute $\widetilde{A},\widetilde{B},\widetilde{C}$ using the formulas~\eqref{eq:tildeABC}\;
apply Algorithm~\ref{algo:crt} to $\widetilde{A},\widetilde{B},\widetilde{C}$ (only the last iterate $\offdiag(\sad{B}_k)$ and $\sad{\bs{v}}= \bs{u} A_k$ are needed)\;
compute $X$ (or $P=h^{-1}Y$, $s=h^{-1}$) and $\Psi$ from~\eqref{cstableinvpairformula}, using the triplet representations derived from Lemma~\ref{subtriplets} to invert $B_{22}$ and $S$\;
compute the triplet representation~\eqref{tripletX2}\;
\caption{Computing a c-stable invariant pair of $P(z)$.} \label{algo:npt}
\end{algorithm}

\begin{remark}
In the case $V=0$, our construction reduces to the method to transform a fluid queue into a QBD introduced by Ramaswami~\cite{ram99}, up to a diagonal scaling. Indeed, the transition matrices $A_0,A_1,A_2$ appearing in~\cite[Equation~4.5]{ram99} satisfy $A_0=K^{-1}\widetilde{C}$, $A_1=I-K^{-1}\widetilde{B}$, $A_2=K^{-1}\widetilde{A}$, where $K$ is the diagonal matrix with entries
  \[
   K_{ii} = \begin{cases}
             d_{ii} & d_{ii}>0,\\
             2\abs{d_{ii}} & d_{ii}<0.
            \end{cases}
  \]
(the case $d_{ii}=0$ is not treated in~\cite{ram99}).
\end{remark}

\subsection{An SDA-like variant}
In the linear case, a popular algorithm for this problem is the~\emph{structured doubling algorithm}~\cite{glx05} (SDA) and its variants~\cite{bmp,wwl}. It is a slightly different iteration, which has a lower computational cost because it uses the block structure in a more effective way. Merging the derivation in~\cite{bmp} with ours, we can obtain a SDA-lookalike variant for second-order problem. The following algorithm indeed reduces to SDA-ss~\cite{bmp} if $n_1=0$.

We start from the matrix polynomial $P(z)$ in the three-blocks form~\eqref{VDQ3blocks}, but this time we apply the discretization map $y=1+hz$ first, and then we modify the location of the infinite eigenvalues. We have $F(y) = Ay^2 + By + C$, with coefficients as in~\eqref{ABC}, that is,

\begin{align*}
A &= \m{h^{-2} V_1\\ & 0 \\ & & 0}, \quad B = \m{h^{-1}D_1+ 2h^{-2}V_1 \\ & h^{-1}D_2\\&& h^{-1}D_3}\\
C &= \m{Q_{11}+ h^{-1}D_1+h^{-2}V_1 & Q_{12} & Q_{13}\\ Q_{21} & Q_{22}+ h^{-1}D_2 & Q_{23}\\  Q_{31} & Q_{32} & Q_{33}+ h^{-1}D_3}.
\end{align*}

We postmultiply these coefficients by the inverse of the M-matrix
\begin{equation} \label{sdalikemmatrix}
M = \m{h^{-2}V_1 &0 &0 \\ 0 & h^{-1}D_2 & 0\\ -Q_{31} & -Q_{32} & -Q_{33}- h^{-1}D_3},	
\end{equation}
an operation which does not change eigenvalues and left invariant pairs, obtaining $\widehat{F}(y) = \widehat{A}y^2 +\widehat{B}y+\widehat{C}$, with
\[
\widehat{A} = \m{I\\ & 0 \\ & & 0}, \quad 
\widehat{B} = \m{B_{11} \\ & I \\-B_{31} & -B_{32} & -B_{33}}, \quad
\widehat{C} = \m{C_{11} & C_{12} & C_{13}\\ C_{21} & C_{22} & C_{23}\\  0 & 0 & -I}, \quad
\]
where the block coefficients are given by
\begin{align*} 
 \m{C_{11} & C_{12} & C_{13} \\
    C_{21} & C_{22} & C_{23}\\
 } & =
 \m{Q_{11}+ h^{-1}D_1+ h^{-2}V_1 & Q_{12} & Q_{13}\\
    Q_{21} & Q_{22}+ h^{-1}D_2 & Q_{23}
 } \times  \\
 & 
\;\;\;\; \m{h^{-2}V_1 &0 &0 \\ 0 & h^{-1}D_2 & 0\\ -Q_{31} & -Q_{32} & -Q_{33}- h^{-1}D_3}^{-1}
\end{align*} 
and
\begin{align*} 
 \m{B_{11} & 0 & 0\\
    B_{31} & B_{32} & B_{33}
 } & =
 \m{h^{-1}D_1+ 2h^{-2}V_1 & 0 & 0\\
    0&0&- h^{-1}D_3
 }  \times \\
 & \;\;\;\;
 \m{h^{-2}V_1 &0 &0 \\ 0 & h^{-1}D_2 & 0\\ -Q_{31} & -Q_{32} & -Q_{33}- h^{-1}D_3}^{-1}.
\end{align*} 
To obtain these blocks with a subtraction-free expression, we can make use of the triplet representation $\left(\offdiag(M),\ones,
\begin{bmatrix}
  h^{-2}\diag(V_1), h^{-1}D_2,-h^{-1}D_3
\end{bmatrix}^{\top}\right)$ for~$M$.

Finally, we postmultiply by
\[
\widehat{S}(y) = \m{I\\&I\\&&yI},
\]
which has the effect of shuffling around some blocks and moving $n_3$ of the infinite eigenvalues to zero; the final result is $\check{F}(y) = \check{A}y^2 + \check{B}y + \check{C}$, with
\begin{equation}  \label{SDAlikeABC}
\check{A}=\m{I & 0 & 0\\0 & 0 & 0 \\ 0 & 0 & B_{33}}, 
\check{B}=\m{B_{11} & 0 & -C_{13}\\0 & I & -C_{23}\\ -B_{31} & -B_{32} & I}, 
\check{C} =\m{C_{11} &C_{12} & 0\\ C_{21} & C_{22} & 0\\0 & 0 & 0}.
\end{equation}
The triple $\check{A}, \check{B}, \check{C}$ has the right signs for us to apply Cyclic Reduction, producing the same solution matrix $R$ as the above approach, since the final location of the eigenvalues is the same. Moreover, some of the pattern in the matrices $\check{A}, \check{B}, \check{C}$ is preserved under CR iterations; namely, at each step $k$, the pattern is
\[
A_k = \m{\ast & \ast & 0\\ \ast & \ast & 0\\ 0 & 0 & 0}, \quad
B_k = \m{\ast & \ast & \ast\\ \ast & I & \ast \\ \ast & \ast & I}, \quad
\sad{B}_k = \m{\ast & \ast & \ast\\ \ast & I & \ast \\ \ast & \ast & I}, \quad
C_k = \m{\ast & 0 & \ast\\ 0 & 0 & 0\\ \ast & 0 & \ast}.
\]

A slightly more efficient version of Cyclic Reduction can be obtained by exploiting the knowledge of these zero and identity blocks. While in the linear case the formulas simplify notably, in our quadratic case it is dubious whether it is worthwhile dealing with the additional complication of these formulas in the implementation, despite the slight computational advantage.

\subsection{Solving the ODE} \label{sec:solvingODE}
The reference~\cite[Sections~1.4, 2.4 and~2.5]{GohLR} contains a complete theory of the relations between invariant pairs and solution of matrix linear differential equations. Let $(X,U)$, with $X\in\mathbb{R}^{\ell\times \ell}$ and $U \in \mathbb{R}^{\ell\times n}$ in the form~\eqref{ivaU}, be a c-stable invariant pair of~\eqref{eq:matpol}. We assume positive recurrence, since otherwise there is no invariant density to compute. Then, the eigenvalues of $X$ coincide with the eigenvalues of $P(z)$ in the open left half-plane, and any solution $\bs{p}(x)$ of~\eqref{ODE} such that $\lim\limits_{x\to\infty} \bs{p}(x)=\bs{0}$ can be written as 
	\begin{align*} 
	\bs{p}(x)=\bs{v}\exp(Xx)U, \quad \mbox{for $\bs{v}\in\mathbb{R}^{1\times \ell}$}.
	\end{align*} 
	
Simple probability considerations show that the invariant measure of the Markov-modulated Brownian motion with coefficients $V,D,Q$ is the sum of a mass $
\begin{bmatrix}
  \bs{0} & \bs{p}_0
\end{bmatrix}$ 
at $x=0$ (where the matrix partitioning is consistent with~\eqref{ivaU}), and the density $\bs{p}(x)=\bs{v}\exp(Xx)U$. If the computed invariant pair satisfies~\eqref{ivaU}, the unknown coefficients $\bs{p}_0$ and $\bs{v}$ can be determined from the condition
\[
\bs{u} = 
\begin{bmatrix}
  \bs{u}_{1} & \bs{u_2}  
\end{bmatrix}
=
\begin{bmatrix}
  0 & \bs{p}_0
\end{bmatrix}
+
\int_{0}^\infty
\bs{p}(x) \mathrm{d}x
=
\begin{bmatrix}
  -\bs{v}X^{-1} & \bs{p}_0 - \bs{v}X^{-1}\Psi
\end{bmatrix}.
\]
Using the relation already derived in~\eqref{vBeq}, we get
\begin{align*}  
\bs{p}_0 = \bs{u}_2 + \bs{v}X^{-1}\Psi = \bs{u}_2 - \bs{u}_1\Psi = \sad{\bs{v}}_2 B_{22}^{-1}.
\end{align*} 
Moreover, $\bs{v}$ satisfies $-\bs{u}_1X=\bs{v}$, hence it follows from~\eqref{tripletX2} that
\[
\bs{v} = \frac{1}{h}\left(\sad{\bs{v}}_1 + \sad{\bs{v}}_2B_{22}^{-1}(C_{21}-B_{21})\right)S^{-1},
\]
the vector that we already have computed when obtaining a triplet representation for $-X^{\top}$.

Hence we have all the quantities that are needed to compute the invariant density $\bs{p}(x)=\bs{v}\exp(Xx)U$.

\section{Componentwise stability} \label{sec:stability}
In this section, we adapt the theory in \cite{NguP15} to prove that the computation of invariant pairs and matrix equation solutions with Cyclic Reduction (in the non-null-recurrent case) is componentwise stable, provided that one uses triplet representations and the GTH trick as described in Algorithm~\ref{algo:crt}.

We define for each $k=0,1,2,\dots$ the 4-tuple of nonnegative matrices and vectors
\begin{equation} \label{Qtuple}
	S_k := (A_k,-\offdiag(B_k),-\offdiag(\sad{B}_k), C_k),
\end{equation}
and we call $\mathcal{F}$ the map such that $S_{k+1}=\mathcal{F}(S_k)$, corresponding to one step of Cyclic Reduction computed with~\eqref{eq:cr}.

When $S_k$ and $\widetilde{S}_k$ are two different 4-tuples in the form~\eqref{Qtuple} and $\alpha$ is a real number, we write for short $\abs{\widetilde{S}_k-S_k} \leq \alpha S_k$ to mean that the relation holds when we replace $S_k$ with each of the matrices and vectors in the 4-tuple, i.e.,
\begin{align*}
	\abs{\widetilde{A}_k-A_k} &\leq \alpha A_k,\\
	\abs{\offdiag(\widetilde{B}_k)-\offdiag(B_k)} &\leq \alpha \abs{\offdiag(B_k)}, \\
	\abs{\offdiag(\widetilde{\sad{B}}_k)-\offdiag(\sad{B}_k)} &\leq \alpha \abs{\offdiag(\sad{B}_k)}, \\
	\abs{\widetilde{C}_k-C_k} &\leq \alpha C_k.
\end{align*}

\subsection{Componentwise perturbation bounds}

We start from assessing the error incurred when starting from inaccurate initial values. We focus on first-order results, and adopt the notation $M \leqdot N$ to mean $M \leq N+O(\varepsilon^2)$.

The key to this result is interpreting the iterates of Cyclic Reduction as the result of a censoring operation. The connection between Cyclic Reduction and censoring is a well-established result (see, for example, \cite[Section~7.3]{blm05}). The following lemma is one of the possible ways to formalize this connection.
\begin{lemma} \label{censorlemma}
Consider the sequences obtained by Cyclic Reduction~\eqref{eq:cr}, and in addition
the sequence $\happy{B}$ defined by $\happy{B}_0 = B_0$ and $\happy{B}_{k+1} =\happy{B}_k - A_k B_k^{-1}C_k$. The matrix
\begin{equation} \label{smallmatrix}
	\begin{bmatrix}
	A_0+I-\sad{B}_k & C_k\\
	A_k & I-\happy{B}_k+C_0
	\end{bmatrix}
\end{equation}
is the result of censoring all blocks apart from the first and last from the $n(2^k+1)\times n(2^k+1)$ matrix
\begin{equation} \label{bigmatrix}
\begin{bmatrix}
  A_0 + I-B_0 & C_0\\
  A_0 & I-B_0 & C_0\\
  & \ddots & \ddots & \ddots\\
  && A_0 & I-B_0 & C_0\\
  &&& A_0 & I-B_0+C_0
\end{bmatrix}
\end{equation}
\end{lemma}
\begin{proof}
We first censor the even-numbered blocks, obtaining
\begin{multline*}
\displaybreak
\begin{bmatrix}
  A_0+I-B_0 \\
  & I-B_0\\
  & & \ddots\\
  & & & I-B_0\\
  & & &  & I-B_0 + C_0\\
\end{bmatrix}  \\
+\begin{bmatrix}
  C_0\\
  A_0 & C_0\\
   & A_0 & \ddots\\
   & & \ddots & C_0\\
   & & & A_0
\end{bmatrix}
\begin{bmatrix}
  B_0\\ & B_0\\ &&\ddots \\ &&& B_0
\end{bmatrix}^{-1}
\begin{bmatrix}
  A_0 & C_0\\
  & A_0 & C_0\\
  && \ddots & \ddots\\
  &&&A_0 & C_0
\end{bmatrix}\\
=  
\begin{bmatrix}
  A_0 + I-\sad{B}_1 & C_1\\
  A_1 & I-B_1 & C_1\\
  &\ddots & \ddots & \ddots\\
  & & A_1& I-B_1 & C_1\\
  & & & A_1 & I-\happy{B}_1 + C_0\\
\end{bmatrix}.
\end{multline*}
We reiterate the same process $k$ times in total, each time censoring the even-numbered blocks in the new matrix; after each step, we obtain a matrix with the same structure, smaller size, and the indices increased by 1.
\end{proof}
\begin{remark}
If the elements in $\diag(B_0)$ are small enough that $I-B_0 \geq 0$, then the matrix in~\eqref{bigmatrix} is stochastic, and so is its censoring~\eqref{smallmatrix}. This gives an alternative proof of the relations 
	\begin{align*}
	(A_k-B_k+C_k)\ones=(A_0-\sad{B}_k+C_k)\ones=(A_k-\happy{B}_k+C_0)\ones = \zeros,
	\end{align*} 
	which appeared in Theorem~\ref{thm:crconv} and Lemma~\ref{lem:secondcrtriplet}.
\end{remark}

Once this lemma is set up, it is simple to prove the following perturbation bound.
\begin{lemma} \label{pertboundlemma}
Let $A,B,C\in \mathbb{R}^{n\times n}$ and $\widetilde{A},\widetilde{B},\widetilde{C}\in \mathbb{R}^{n\times n}$
be two different triples of matrices satisfying A4, A5, A6, such that
\begin{subequations}
\begin{align}
	\abs{\widetilde{A}-A} &\leqdot \varepsilon A, \\
	\abs{\offdiag(\widetilde{B})-\offdiag(B)} &\leqdot \varepsilon \abs{\offdiag(B)}, \\
	\abs{\widetilde{C}-C} &\leqdot \varepsilon C.
\end{align}
\end{subequations}
Let $S_k$ and $\widetilde{S}_k$ be the 4-tuples resulting from applying $k$ steps of Cyclic Reduction~\eqref{eq:cr} starting from $A,B,C$ and $\widetilde{A},\widetilde{B},\widetilde{C}$, respectively. Then,
\[
\abs{\widetilde{S}_k-S_k}\leqdot n2^k \varepsilon S_k.
\]
\end{lemma}
\begin{proof}
Up to a common scaling factor (which does not alter the statement of the theorem), we can assume that $I-B\geq 0$. Then, the matrices in~\eqref{smallmatrix} and~\eqref{bigmatrix} are stochastic, and we can apply~\cite[Lemma~7.3]{NguP15} to this censoring operation. 

In detail, we call $P$ the matrix in~\eqref{bigmatrix}, and $\widetilde{P}$ its equivalent built starting with the initial values with a tilde. We have for $i\neq j$
\[
\abs{(\widetilde{A}_0+I-\widetilde{B}_0)_{ij} - (A_0+I-B_0)_{ij}} \leqdot \abs{(\widetilde{A}_0-A_0)_{ij}} + \abs{(\widetilde{B}_0-B_0)_{ij}} \leqdot \varepsilon (A_0+I-B_0)_{ij},
\]
and similarly for all other entries, so $\abs{\offdiag(\widetilde{P})-\offdiag(P)}\leq \varepsilon\offdiag(P)$.
Thus, the first part of~\cite[Lemma~7.3]{NguP15} holds with $m = n(2^k-1)\leq n2^k$. This proves that
\begin{align*}
	\abs{\widetilde{A}_k-A_k} &\leqdot n2^k \varepsilon A_k,\\
	\abs{\widetilde{C}_k-C_k} &\leqdot n2^k \varepsilon C_k,\\
	\intertext{and}
	\abs{\widetilde{D}_k-D_k} &\leqdot n2^k\varepsilon D_k,\\
	\abs{\widetilde{E}_k-E_k} &\leqdot n2^k\varepsilon E_k,
\end{align*}
where
\begin{align*}
D_k & = B_0-\sad{B}_k = \sum_{j=0}^{k-1} C_jB_j^{-1}A_j, \quad E_k = B_0-\happy{B}_k = \sum_{j=0}^{k-1} A_jB_j^{-1}C_j, 
\end{align*} 
	and equivalent definitions with the tilde symbols. The bounds
\begin{align*}
\abs{\offdiag(\widetilde{B}_k)-\offdiag(B_k)} &\leqdot n2^k\varepsilon \abs{\offdiag(B_k)},\\
\abs{\offdiag(\widetilde{\sad{B}}_k)-\offdiag(\sad{B}_k)} &\leqdot n2^k\varepsilon \abs{\offdiag(\sad{B}_k)}
\end{align*}
follow by noting that $B_k=B_0-D_k-E_k$, $\sad{B}_k=B_0-D_k$ and using \cite[Lemma~7.2 (i)]{NguP15}.
\end{proof}

\subsection{Stability of a CR step}
Our next point is investigating the stability of a step of Cyclic Reduction when performed in machine arithmetic. We rely once again on the lemmas on basic operations in~\cite[Section~7]{NguP15}, and we hide in $M\leqdot N$ terms which are second-order in $\mp$.

\begin{lemma} \label{onecrstep}
Let the $4$-tuple $S_k = (A_k,-\offdiag(B_k),-\offdiag(\sad{B}_k),C_k)$ be exactly-represented machine numbers. We denote by $S_{k+1} =\mathcal{F}(S_k)$ the result of performing one step of Cyclic Reduction on them, and by $\widetilde{S}_{k+1} = \widetilde{\mathcal{F}}(S_k)$ the result of performing one step of Cyclic Reduction computed in inexact machine arithmetic, starting from the same matrices.

Then,
\[
\abs{\widetilde{S}_{k+1}-S_{k+1}}\leqdot (\psi(n)+n+2)\mp S_{k+1},
\]
where $n$ is the size of the involved matrices and $\psi(n)=\frac{2}{3}(2n+5)(n+2)(n+3)$ is the accuracy bound for the solution of a linear system with the GTH algorithm (as in~\cite[Theorem~4.1]{NguP15}).
\end{lemma}
\begin{proof}
We use, with a slight abuse of notation, the notation $c(X)$ to denote the computed approximation of a quantity $X$ along one step of the algorithm (even though it is not, strictly speaking, a function of $X$ only).

Using \cite[Lemma~7.9]{NguP15} with $a=b=0$, we obtain that the computed value $c(B_k^{-1}A_k)$ of $B_k^{-1}A_k$ satisfies
\[
\abs{c(B_k^{-1}A_k)-B_k^{-1}A_k} \leqdot \psi(n)\mp B_k^{-1}A_k.
\]
Hence the computed values of $A_{k+1}=A_kB_k^{-1}A_k$ and $C_{k}B_k^{-1}A_k$ satisfy (by \cite[Lemma~7.8]{NguP15})
\begin{align*}
\abs{c(A_{k+1})-A_{k+1}} &\leqdot (\psi(n)+n)\mp A_{k+1},\\
\abs{c(C_{k}B_k^{-1}A_k)-C_{k}B_k^{-1}A_k} &\leqdot (\psi(n)+n)\mp C_{k}B_k^{-1}A_k,\\
\end{align*}
Analogously we have
\begin{align*}
\abs{c(C_{k+1})-C_{k+1}} &\leqdot (\psi(n)+n)\mp C_{k+1},\\
\abs{c(A_{k}B_k^{-1}C_k)-A_{k}B_k^{-1}C_k} &\leqdot (\psi(n)+n)\mp A_{k}B_k^{-1}C_k,\\
\end{align*}
Using again \cite[Lemma~7.8]{NguP15} for the additions, we have then
\begin{align*}
\abs{c(\offdiag(B_{k+1}))-\offdiag(B_{k+1})} &\leqdot (\psi(n)+n+2)\mp \abs{\offdiag(B_{k+1})},\\
\abs{c(\offdiag(\sad{B}_{k+1}))-\offdiag(\sad{B}_{k+1})} &\leqdot (\psi(n)+n+1)\mp \abs{\offdiag(\sad{B}_{k+1})}.
\end{align*}
\end{proof}

\subsection{Stability of multiple CR steps}
We can now address multiple steps of Cyclic Reduction. The proof here follows~\cite[Theorem~7.12]{NguP15}.
\begin{lemma}
Let $A,B,C\in\mathbb{R}^{n\times n}$ be three matrices satisfying Assumptions A4, A5, A6, and such that $A,C$ and $\offdiag(B)$ are exactly-represented machine numbers. Denote by $S_k=\mathcal{F}^k(S_0)$ the result of performing $k$ steps of Cyclic Reduction starting from $S_0=(A,-\offdiag(B),-\offdiag(B),-\offdiag(B),C)$, and by $\widetilde{S}_k=\widetilde{\mathcal{F}}^k(S_0)$ the result of $k$ steps of Cyclic Reduction performed in inexact machine arithmetic. Then,
\begin{equation} \label{crstabequation}
	\abs{\widetilde{S}_k-S_k} \leq n2^k(\psi(n)+n+2)\mp S_k.
\end{equation}
\end{lemma}
\begin{proof}
We prove the result by induction on $k$; the base case ($k=1$) is Lemma~\ref{onecrstep}.

The following manipulation is a formal version of the statement that when considering first-order error bounds we can add up the local errors at the different steps of the algorithm. Consider the telescopic sum
\begin{equation} \label{telesum}
\abs{\widetilde{S}_k-S_k} \leq \sum_{h=1}^k \abs{\mathcal{F}^{h-1}\widetilde{\mathcal{F}}(\widetilde{S}_{k-h}) - \mathcal{F}^{h-1}\mathcal{F}(\widetilde{S}_{k-h})}.
\end{equation} 
By Lemma~\ref{onecrstep}, we have
\[
\abs{\widetilde{\mathcal{F}}(\widetilde{S}_{k-h})-\mathcal{F}(\widetilde{S}_{k-h})} \leqdot (\psi(n)+n+2)\mp \mathcal{F}(\widetilde{S}_{k-h}).
\]
Then by Lemma~\ref{censorlemma} used with $\varepsilon=(\psi(n)+n+2)\mp$, we get
\begin{align*}
\abs{\mathcal{F}^{h-1}\widetilde{\mathcal{F}}(\widetilde{S}_{k-h}) - \mathcal{F}^{h-1}\mathcal{F}(\widetilde{S}_{k-h})} &\leqdot n2^{h-1}(\psi(n)+n+2)\mp \mathcal{F}^{h-1}\mathcal{F}(\widetilde{S}_{k-h}) \\
&\leqdot  n2^{h-1}(\psi(n)+n+2)\mp \mathcal{F}^h(S_{k-h}) \\
&= n2^{h-1}(\psi(n)+n+2)\mp S_{k}.
\end{align*}
Passing from the first to the second row we have replaced $\widetilde{S}_{k-h}$ with $S_{k-h}$; this is possible because they differ by a term of order $O(\mp)$ by inductive hypothesis.

Insert this inequality into~\eqref{telesum} to get
\[
\abs{\widetilde{S}_k-S_k} \leqdot \sum_{h=1}^k n2^{h-1}(\psi(n)+n+2)\mp S_{k} < n2^k(\psi(n)+n+2)\mp S_{k}.
\]
\end{proof}

\subsection{Putting everything together}
The previous sections shows that the CR iteration (Algorithm~1) is componentwise stable. The computation of its initial values starting from $V,D,Q$ can be performed with~\eqref{ABC}, \eqref{eq:tildeABC}, or~\eqref{SDAlikeABC}; in all three cases, if~\eqref{eq:hconstraints} holds for each $i\in E_1\cup E_2$, then we obtain an approximation $\widetilde{S}_0$ of the CR initial values satisfying $\abs{\widetilde{S}_0-S_0}\leq \alpha \mp S_0$ for a moderate multiple $\alpha$ of the machine precision, and thus by Lemma~\ref{pertboundlemma} the computed iterates are also componentwise accurate.

Once a sufficient number $k$ of steps is performed to achieve convergence, we compute the invariant pair $(X,U)$ as described in Section~\ref{sec:deflatingzeros}. The computed iterates satisfy~\eqref{crstabequation}, and similarly the computed approximation $\widetilde{\sad{\bs{v}}}_k$ of $\sad{\bs{v}}=\sad{\bs{u}}A_k$ satisfies
\[
 \abs{\widetilde{\sad{\bs{v}}}_k-\sad{\bs{v}}_k} = \abs{\bs{u}\widetilde{A}_k- \bs{u}A_k} \leqdot n2^k(\psi(n)+n+2)\mp \sad{\bs{v}}_k.
\]

The rest of the computation only involves subtraction-free formulas: we have described in Sections~\ref{sec:deflatingzeros}, \ref{sec:tripletX2}, and~\ref{sec:solvingODE} how to get from the matrix $R$ computed by CR (and $\sad{B}_{\infty}$ and $\sad{\bs{v}}$) to the stable invariant pair of $P(z)$, its triplet representation, and the quantities needed to compute $\bs{p}(x)$.

\section{Numerical experiments} \label{sec:numerical}
We compare the following methods.
\begin{description}	
	\item[KK] The algorithm in~\cite{kk95}, based on explicit computation of eigenvalues and eigenvectors of a linearizing matrix which is obtained (essentially) by deflating the infinite eigenvalues from the linearizing matrix polynomial
\begin{equation} \label{linearization}
\mathcal{A} - z\mathcal{E} = 
	\begin{bmatrix}
	  D & -T^\top\\
	  I_n & 0
	\end{bmatrix}
	- z
	\begin{bmatrix}
	  V & 0\\
	  0 & I_n
	\end{bmatrix}.
\end{equation}
	The main drawback of this method, is that by computing explicitly an eigenvalue decomposition we expect error amplification by the condition number of the eigenvector matrix.
	\item[AS] The algorithm in~\cite{AgaS}, based on computing the sign function of the pencil~\eqref{linearization} using the Newton-like iteration $\mathcal{A}_{k+1} = \frac{1}{2}\left( \mathcal{A}_k + \mathcal{E}\mathcal{A}_k^{-1}\mathcal{E}\right)$, and using it to separate the infinite, stable and unstable eigenvalues into different blocks. In principle, this algorithm goes in the right direction to get better numerical properties; in practice, unfortunately, our implementation of this algorithm was affected negatively by convergence issues in this iteration. In particular, it seems that the pencil $\mathcal{A}_k-z\mathcal{E}$ converges to a singular pencil whenever $\diag(V)$ has zero entries, so the inversion $\mathcal{A}_k^{-1}$ becomes increasingly ill-conditioned. This complicates the choice of a stopping criterion.
	\item[QZ] An algorithm similar to AS, but in which the stable subspace is computed using a permuted QZ decomposition~\cite{Kag92} of~\eqref{linearization} (MATLAB's \texttt{ordqz}). While we could not find an explicit reference in the applied probability literature for the use of this method in the context of Markov-modulated Brownian motion, it is the method of choice for problems of this kind in the numerical linear algebra community \cite{BetK}. The QZ decomposition is normwise backward stable, so we expect excellent normwise stability properties.
	\item[LN] The algorithm in~\cite{LatN}, which is based on Cyclic Reduction without the use of triplet representations, or of any particular method to preserve positiveness. The discretizing transformation is the Cayley transform with $h=1$, $y = (z+1)/(z-1)$. This algorithm can solve only problems with $\diag(V)>0$.
	\item[NP] Algorithm~\eqref{algo:npt} as described here, from which we expect componentwise accuracy.
\end{description}
We apply these algorithms to several test problems.
\begin{description}
	\item[NP15] A modification of \cite[Example~5.1]{NguP15}, a problem in which there is an imbalance of several orders of magnitude between the components of the solution. We take $T$ and $D$ as in that problem, and add a Brownian motion component with $V=I$.
	\item[NP15s] The same problem as NP, but with $V(n,n)=0$, to obtain a problem with singular $V$.
	\item[rand($n$)] Random-generated problems of different sizes $n = 8,20,50$. The matrices are generated with the MATLAB commands
	\begin{verbatim}
V = diag(abs(randn(n, 1));
D = diag(randn(n, 1));
T = abs(randn(n)); T = T - diag(T * ones(n, 1));
	\end{verbatim}
	\item[rand($n$)s] The same as rand($n$), but with a matrix $V$ containing four zero diagonal entries: \texttt{V = blkdiag(diag(abs(randn(n-4, 1))), zeros(4))}.
	\item[imb($n$),imb($n$)s] Defined as rand($n$) and rand($n$)s, but all the calls of the form \texttt{randn(h,k)} are replaced by a different procedure that generates numbers spanning different orders of magnitude: \texttt{randn(h, k) .* exp(5 * randn(h, k))}.
\end{description}
To improve reproducibility without generating the same numbers repeatedly, we have reset the random number seed once before the complete set of experiments.

As a first error measure, we have considered the residual in the Euclidean norm
\begin{equation} \label{relresidual}
	\frac{\norm{X^2UV - XUD + UQ}}{\norm{U}(\norm{V}+\norm{D}+\norm{Q})}
\end{equation}
of the left stable invariant pair $(X,U)$ as produced by the algorithms. The values of this residual are in Table~\ref{table:residuals}.

Moreover, we have normalized each invariant pair to be in the form~\eqref{ivaU} with a similarity transformation~\eqref{eq:qtrans}, and checked the forward errors
\begin{equation} \label{forwarderror}
\frac{\norm{X-X_{\mathrm{exact}}}}{\norm{X_{\mathrm{exact}}}}, \quad \frac{\norm{\Psi-\Psi_{\mathrm{exact}}}}{\norm{\Psi_{\mathrm{exact}}}},
\end{equation}
where the reference values $X_\mathrm{exact}$ and $\Psi_\mathrm{exact}$ are computed applying method KK with higher precision arithmetic (32 digits, using Matlab's \texttt{vpa} command). The results are in Tables~\eqref{table:Xrel} and~\eqref{table:Psirel}. Note that in the problems with $V>0$, we have $E_3=\varnothing$, and hence the matrix $\Psi$ is empty and computing the error does not make sense. For this reason, Table~\eqref{table:Psirel} does not contain all the experiments.
\begin{table}                                          
\centering                                                 
\begin{tabular}{rccccc}
\toprule                              
 Problem & KK & AS & LN & QZ & NP \\                               
\midrule
NP15 & 1.5e-15 & 5.7e-08 & 3.7e-15 & 9.9e-16 & 3.8e-16 \\    
NP15s & 5.0e-16 & 5.0e-08 & - & 7.7e-16 & 2.3e-16 \\         
rand8 & 1.5e-15 & 9.4e-16 & 3.6e-15 & 9.9e-16 & 1.1e-15 \\ 
rand8s & 6.6e-15 & 2.8e-11 & - & 5.5e-15 & 2.6e-15 \\      
rand20 & 1.9e-15 & 5.6e-14 & 1.8e-14 & 1.8e-15 & 7.3e-16 \\
rand20s & 2.8e-15 & 5.8e-12 & - & 2.0e-14 & 1.3e-14 \\     
rand50 & 2.3e-15 & 3.0e-14 & 3.4e-13 & 1.5e-14 & 5.9e-15 \\
rand50s & 7.1e-14 & 1.2e-08 & - & 1.3e-14 & 1.7e-14 \\     
imb8 & 1.2e-08 & 8.6e-05 & 4.5e+04 & 4.2e-11 & 7.4e-09 \\  
imb8s & 1.5e-13 & 4.2e-10 & - & 2.3e-14 & 2.3e-13 \\       
imb20 & 2.2e-15 & 4.8e-06 & 9.7e-01 & 2.7e-14 & 4.9e-13 \\ 
imb20s & 1.2e-09 & 8.0e+05 & - & 2.6e-13 & 1.9e-13 \\      
imb50 & 8.7e-14 & 7.5e-06 & 3.3e+01 & 3.3e-13 & 1.3e-10 \\ 
imb50s & 3.1e-04 & 3.2e+11 & - & 2.5e-05 & 2.0e-08 \\      
\bottomrule
\end{tabular}                                              
\caption{Relative residual~\eqref{relresidual}.}                                       
\label{table:residuals}
\end{table}                                                
\begin{table}
\centering                                                 
\begin{tabular}{rccccc}                                    
\toprule                              
 Problem & KK & AS & LN & QZ & NP \\                               
\midrule
NP15 & 2.7e-12 & 2.5e-07 & 2.9e-13 & 1.8e-12 & 1.7e-16 \\    
NP15s & 1.3e-12 & 2.2e-07 & - & 6.2e-13 & 1.8e-16 \\         
rand8 & 2.8e-15 & 1.5e-15 & 1.6e-15 & 2.4e-15 & 2.7e-16 \\ 
rand8s & 2.9e-15 & 1.8e-13 & - & 2.3e-15 & 3.1e-16 \\      
rand20 & 4.4e-15 & 9.6e-14 & 5.6e-15 & 4.8e-15 & 3.0e-16 \\
rand20s & 3.2e-15 & 3.0e-12 & - & 4.1e-14 & 1.1e-15 \\     
rand50 & 5.9e-15 & 4.0e-14 & 4.0e-14 & 5.6e-14 & 6.9e-16 \\
rand50s & 5.6e-14 & 1.2e-10 & - & 3.5e-14 & 5.2e-16 \\     
imb8 & 9.7e-12 & 1.9e-09 & 1.1e+00 & 7.1e-13 & 9.0e-13 \\  
imb8s & 2.6e-14 & 1.3e-08 & - & 1.3e-12 & 1.1e-15 \\       
imb20 & 4.6e-11 & 2.1e-07 & 3.2e-04 & 1.1e-09 & 9.1e-12 \\ 
imb20s & 4.4e-12 & 6.9e-06 & - & 5.9e-12 & 4.0e-13 \\      
imb50 & 2.0e-10 & 9.8e-06 & 7.2e-01 & 1.0e-08 & 8.3e-10 \\ 
imb50s & 2.0e-10 & 3.3e-05 & - & 1.0e+00 & 2.6e-13 \\      
\bottomrule
\end{tabular}                                              
\caption{Forward error~\eqref{forwarderror} on $X$.}                                  
\label{table:Xrel}
\end{table}                                                
\begin{table}                                     
\centering                                            
\begin{tabular}{rccccc}                                    
\toprule                              
 Problem & KK & AS & LN & QZ & NP \\                               
\midrule
NP15s & 2.3e-15 & 1.8e-11 & - & 2.8e-15 & 1.3e-16 \\    
rand8s & 1.2e-14 & 3.7e-13 & - & 2.4e-15 & 2.5e-15 \\ 
rand20s & 7.1e-15 & 7.7e-11 & - & 6.7e-14 & 2.1e-15 \\
rand50s & 3.4e-14 & 3.5e-09 & - & 5.3e-14 & 4.7e-16 \\
imb8s & 8.3e-15 & 5.2e-09 & - & 1.1e-11 & 5.2e-15 \\  
imb20s & 1.4e-10 & 1.9e-08 & - & 2.8e-11 & 4.0e-11 \\ 
imb50s & 6.9e-11 & 9.0e-09 & - & 1.0e-04 & 6.1e-08 \\ 
\bottomrule
\end{tabular}                                         
\caption{Forward error~\eqref{forwarderror} on $\Psi$.
}
\label{table:Psirel}
\end{table}

As one can see, the results obtained by the new algorithm are very satisfying, especially in terms of accuracy of the computed $X$ and $\Psi$ (which are often the quantities of interest in view of their physical interpretation). The relative residual of the computed invariant pair, however, is sometimes slightly higher than the one obtained with the QZ method.

\section{Conclusions} \label{sec:conclusions}
We have described a subtraction-free algorithm to compute the quantities needed to determine the steady-state behavior of Markov-modulated Brownian motion models in a componentwise accurate fashion. The algorithm extends the one described in \cite{NguP15} for the linear case $V=0$, and is based on a componentwise accurate variant of Cyclic Reduction. A componentwise error analysis of this CR algorithm is provided. Our analysis highlights the role of the spectral transformation which converts continuous-time to discrete-time stability. Another interesting result is the use of a transformation related to index reduction for differential-algebraic equation and to the  shift technique in this novel context. 

\bibliographystyle{alpha}
\bibliography{NP2013}

\end{document}